\documentclass[11pt,reqno]{amsart}
\usepackage[T1]{fontenc}
\usepackage{lmodern}
\usepackage{a4wide}
\usepackage{microtype}
\usepackage{geometry}
\usepackage{comment}
\usepackage{xcolor}
\usepackage{mathtools,amssymb,amsthm,mathrsfs}
\usepackage{enumitem}
\usepackage[colorlinks=true,linkcolor=blue!45!black,citecolor=blue!45!black, urlcolor=blue!45!black]{hyperref}
\usepackage[nameinlink,noabbrev]{cleveref}
\allowdisplaybreaks

\expandafter\providecommand\csname subjclassname@2020\endcsname
  {\textup{2020} Mathematics Subject Classification}

\newtheorem{theorem}{Theorem}[section]
\newtheorem{proposition}[theorem]{Proposition}
\newtheorem{lemma}[theorem]{Lemma}
\newtheorem{corollary}[theorem]{Corollary}

\theoremstyle{definition}

\theoremstyle{remark}
\newtheorem{remark}[theorem]{Remark}
\numberwithin{equation}{section}

\newcommand{\R}{\mathbb{R}}
\newcommand{\Hh}{\mathcal{H}}
\newcommand{\cA}{\mathcal{A}}
\newcommand{\cE}{\mathcal{E}}
\newcommand{\one}{\mathbf{1}}

\newcommand{\dist}{\operatorname{dist}}
\newcommand{\Jint}{\operatorname{Int}}
\newcommand{\Hmone}[1]{\lVert #1\rVert_{\dot H^{-1}_D}}
\newcommand{\p}{\partial}
\newcommand{\G}{\mathcal G}
\newcommand{\M}{\mathcal M}
\newcommand{\dd}{\,\mathrm d\mathbf{x}}
\newcommand{\ds}{\,\mathrm ds}
\newcommand{\eps}{\varepsilon}

\title[Energy-maximizing vortex patches]
{Boundary regularity, finiteness and stability of energy-maximizing vortex patches}

\author{Daomin Cao}
 \address{State Key Laboratory of Mathematical Sciences, Academy of Mathematics and Systems
Science, Chinese Academy of Sciences, Beijing 100190, P.R. China, and University of Chinese Academy of Sciences, Beijing 100049,  P.R. China}
\email{dmcao@amt.ac.cn}

\author{Guodong Wang}
\address{School of Mathematical Sciences, Dalian University of Technology, Dalian 116024, P.R. China}
\email{gdw@dlut.edu.cn}

\author{Yuchen Wang}
\address{School of Mathematics, Tianjin Normal University \\ Tianjin, China}
\email{wangyuchen@mail.nankai.edu.cn}


\begin{document}

\begin{abstract}
We establish boundary regularity for energy-maximizing vortex patches of the 2-D incompressible Euler equations in a bounded Jordan domain. These results resolve Turkington’s free-boundary regularity question throughout the full range of admissible parameters.  For bounded convex domains, we further prove that the set of energy maximizers is finite at each prescribed patch area. When the domain is smooth, the number of maximizers is uniformly bounded over all admissible patch areas, and every maximizer is Lyapunov stable in $L^p$ for any $1< p<\infty$.
\end{abstract}

\maketitle

\section{Introduction}\label{sec:introduction}

Consider the two-dimensional incompressible Euler equations in a bounded Jordan domain $D\subset\R^2$:
\begin{equation}\label{eq:euler-velocity}
 \begin{cases}
 \partial_t\mathbf v+(\mathbf v\cdot\nabla)\mathbf v=-\nabla p
       &\text{in }D,\\
 \nabla\cdot\mathbf v=0&\text{in }D,
 \end{cases}
\end{equation}
where $\mathbf v=(v_1,v_2)^{\mathsf T}$ is the fluid velocity and $p$ is the pressure. For $\mathbf v(t)\in L^2(D;\R^2)$, incompressibility and the
slip boundary condition are imposed in the weak form
\begin{equation}\label{eq:weak-impermeability}
 \int_D\mathbf v(t,\mathbf x)\cdot\nabla\varphi(\mathbf x)
 \,\mathrm d\mathbf x=0
 \qquad  \forall \varphi\in H^1(D),  t \geq 0. 
\end{equation}
It is not hard to verify that this equation recovers the usual slip boundary condition $\mathbf v\cdot\mathbf n=0$ for a smooth domain, where $\mathbf n$ is the outward unit normal.

Under mild regularity assumptions, it is convenient to work with the vorticity formulation of the Euler equations:
\begin{equation}\label{eq:euler-vorticity}
	\partial_t\omega+\mathbf v\cdot\nabla\omega=0,
	\qquad \mathbf v=\nabla^\perp\mathcal G\omega,
\end{equation}
where $\omega:=\partial_1v_2-\partial_2v_1$ is the vorticity, $\nabla^\perp=(\partial_2,-\partial_1)$ and
$\mathcal G=(-\Delta_D)^{-1}$ is the Green operator for $-\Delta$ with the homogeneous Dirichlet boundary condition. We note that 
\[
 \psi=\mathcal G\omega\in H_0^1(D),\; \forall \omega \in H^{-1}(D)
\]
and $\mathbf{v}=\nabla^\perp \psi$ satisfies \eqref{eq:weak-impermeability}.

Yudovich's theorem \cite{Yudovich1963} asserts that $L^\infty$ vorticity defines an Euler flow in smooth bounded domain globally in time.
Thus the vorticity is transported by the Euler flow. Specifically, if the vorticity distribution takes a patch form initially, the solution retains the patch form for all time. Consequently, the vortex patch dynamics is reduced to the boundary dynamics of vortical domains. The global regularity of vortex patches in the plane was first proved by Chemin in \cite{Chemin1993} provided the initial patch is of class $C^{1,\alpha}$, $0<\alpha<1$, see also \cite{BertozziConstantin,Serfati1994} for another elementary proof. Global regularity and small-scale formation in smooth bounded domains are studied in \cite{KiselevLi}, whereas the ill-posedness of 
$C^2$-vortex patch is established in \cite{KiselevLuo}.
 
A vortex patch is called steady when its vorticity is independent of
time, or equivalently
\[
 \int_D\omega\,\nabla^\perp\mathcal G\omega\cdot\nabla\varphi
 \,\mathrm d\mathbf x=0,
 \qquad \varphi\in C_c^\infty(D).
\]
Following the Kelvin-Benjamin variational principle that steady flows arise as
critical points of kinetic energy over rearrangements of vorticity,
Turkington \cite{TurkingtonI} constructed steady vortex patches by
maximizing the kinetic energy
  \begin{equation}\label{eq:kinetic-energy}
 	\cE_{\mathrm{kin}}(\omega)
 	=\frac12\int_D|\mathbf v|^2\,\mathrm d\mathbf{x}	=\frac12\int_D\omega\mathcal G\omega\,\mathrm d\mathbf{x}
 \end{equation}
 on the admissible class
 \[
 \cA_\lambda=\left\{\omega\in L^\infty(D):
 0\leq\omega\leq\lambda,\quad
 \int_D\omega\,\mathrm d\mathbf{x}=1\right\}.
 \]
Here $\lambda>|D|^{-1}$ is a prescribed vorticity bound, while the circulation is normalized to one. We note that $\cA_{\lambda}$ consists of $\frac{1}{|D|} \one_D$ only when $\lambda = |D|^{-1}$.

In the rearrangement point of view
\cite{Burton1989,Burton2005}, the class of patches
\[
 \mathscr R_\lambda
 :=\{\lambda\one_B:B\subset D\text{ measurable},\ |B|=\lambda^{-1}\}
\]
is a rearrangement class, and $\cA_\lambda$ is its relaxed admissible
class; more precisely, $\cA_\lambda$ is the weak $L^2(D)$ closure of $\mathscr R_\lambda$.
In the smooth-domain setting, Turkington \cite{TurkingtonI}
proved that the maximum is attained and that every maximizer is a steady vortex patch of the form
\[
 \cE_{\mathrm{kin}}(\omega_\lambda)
 =\max_{\omega\in\cA_\lambda}\cE_{\mathrm{kin}}(\omega),
 \qquad \omega_\lambda=\lambda\one_{A_\lambda},
\]
where
\[
 A_\lambda=\{\mathcal G\omega_\lambda>\mu_\lambda\}\subset D,
 \qquad \mu_\lambda>0,\qquad |A_\lambda|=\lambda^{-1}.
\]
The threshold $\mu_\lambda$ is the Lagrange multiplier associated with
the circulation constraint, and we refer to these maximizers as
\emph{energy-maximizing vortex patches}.
For sufficiently large $\lambda>0$, Turkington proved that
$\partial A_\lambda$ is a single simple closed $C^1$ curve and that the
rescaled patches converge to the unit disk
\cite{TurkingtonI}. For arbitrary
$\lambda>|D|^{-1}$, he conjectured that the set of irregular free
boundary points $\Sigma_\lambda:=\{x\in\partial A_\lambda: \nabla\mathcal G\omega_\lambda(x)=0\}$
would satisfy
\[
 \#\Sigma_\lambda<\infty
 \quad\text{or}\quad \mathcal H^1(\Sigma_\lambda)=0.
\]

Our aim is to establish boundary regularity for every energy-maximizing
vortex patch in $\cA_\lambda$. The existence and threshold argument are given in Proposition \ref{prop:turkington}. We denote the bounded component of $\R^2\setminus\Gamma$ by $\Jint(\Gamma)$ when $\Gamma$ is a Jordan curve, and use $\bigsqcup$ to address a pairwise disjoint union. The first result of this paper is given as follow:
\begin{theorem}\label[theorem]{thm:main}
	Let $D\subset\R^2$ be a bounded Jordan domain and $\lambda>|D|^{-1}$.
	The energy functional attains its maximum on
	$\cA_\lambda$, and every maximizer has the representation
	\[
	\omega_\lambda=\lambda\one_{A_\lambda}\quad\text{a.e.},\qquad
	u_\lambda=\mathcal G\one_{A_\lambda}-q_\lambda,
	\qquad q_\lambda>0,
	\]
	where $A_\lambda=\{u_\lambda>0\}$, $\overline{A_\lambda}\subset D$, and
	\[
	-\Delta u_\lambda=\one_{A_\lambda}\quad\text{in }D,
	\qquad u_\lambda=-q_\lambda\quad\text{on }\partial D.
	\]
	The zero level-set of $u_\lambda$ coincides with $\p A_\lambda$. There are finitely many pairwise disjoint real-analytic Jordan curves
	$\Gamma_1,\ldots,\Gamma_N$, $N\geq1$, such that
	\[
	A_\lambda=\bigsqcup_{j=1}^N\Jint(\Gamma_j),\qquad
	\partial A_\lambda=\bigsqcup_{j=1}^N\Gamma_j,
	\qquad \min_{\partial A_\lambda}|\nabla u_\lambda|>0.
	\]
	For each $0<\alpha<1$, there are $C,r_*>0$, depending on the fixed
	patch and $\alpha$, such that
	\begin{equation}\label{eq:main-density}
		\sup_{\mathbf p\in\partial A_\lambda}
		\left|\frac{|A_\lambda\cap B_r(\mathbf p)|}{\pi r^2}-\frac12\right|
		\leq Cr^\alpha\qquad(0<r<r_*).
	\end{equation}
\end{theorem}
\begin{remark}\label[remark]{rem:fluid-boundary}
A $C^0$ boundary is understood in the local graph sense: after a change of coordinates, for a neighborhood $U$ and a continuous
function $h$,
\[
 D\cap U=\{(x_1,x_2)\in U:x_2>h(x_1)\},\qquad
 \partial D\cap U=\{(x_1,x_2)\in U:x_2=h(x_1)\}.
\]

\end{remark}

The super level-set representation of  energy-maximizing vortex patches reduces their boundary regularity to the unstable
obstacle equation
\begin{equation}\label{eq:unstable-intro}
	-\Delta u=\one_{\{u>0\}},
\end{equation}
whose singular solutions and variational regularity are studied in
\cite{AnderssonWeiss,AnderssonShahgholianWeissInvent,
AnderssonShahgholianWeissCMP,MonneauWeiss}.
With the normalized energy
$\cE(\rho)=\tfrac12\langle\rho,\mathcal G\rho\rangle$ and admissible class
\begin{equation}\label{eq:mass-class}
 \mathcal A_m=\left\{\rho\in L^\infty(D):0\leq\rho\leq1, \; \int_D\rho=m\right\}.
\end{equation}
By the change of normalization 
\[
 m=\lambda^{-1},\qquad \omega=\lambda\rho,\qquad
 \cE_{\mathrm{kin}}(\omega)=\lambda^2\cE(\rho),
\] 
the normalized stream function $u=\lambda^{-1}(\mathcal G\omega_\lambda-\mu_\lambda)$ satisfies equation\eqref{eq:unstable-intro}. For the boundary regularity of uniformly
rotating vortex patches, one has a sign-changing elliptic free-boundary problem, see for example
\cite{HmidiMateuVerdera,CastroCordobaGomezSerrano,WangZhangZhou}.
\begin{remark}
At the endpoint $\lambda=|D|^{-1}$, the admissible class consists only of
$\lambda\one_D$, whose boundary inherits the
regularity of $D$. For every nontrivial bound $\lambda>|D|^{-1}$,
Theorem \ref{thm:main} gives $\Sigma_\lambda=\varnothing$ and hence completely resolves Turkington's free-boundary regularity question. The theorem concerns each connected component
and does not assert that the patch is connected; see Remark \ref{rem:connectedness}.
\end{remark}

Our approach follows the weighted trace and annular arguments developed by Chanillo, Kenig,
and To \cite{ChanilloKenigTo} for the composite membrane problem. For the nonlocal kinetic energy,
we derive the area-constrained second-variation inequality \eqref{eq:shape} from the exact quadratic energy expansion. The variations act on arbitrary finite collections of regular boundary
arcs under a single total-area constraint, retaining the Green interactions between distinct pieces. Strong convergence of the
rescaled-vorticity layers in $H^{-1}(D)$ justifies passage to the
quadratic term. The resulting weighted trace inequality yields annular
summability, while a torsion identity gives the scale-independent lower
bound $4\pi$ for the weighted boundary integral of each closed regular component. Together with disk comparison, these estimates exclude singular points on the free boundary for every $m\in(0,|D|)$. Recently, 
Ferreri, Mazari-Fouquer, and Prunier \cite{FerreriMazariPrunier} establish
regularity for volume-constrained optimal control problems with bilinear state equations of the form
$-\Delta\Theta_m=m\Theta_m+B(x,\Theta_m)$, while the present problem has the fixed Poisson state equation $-\Delta\psi=\rho$ and the quadratic
functional $\cE(\rho)=\frac12\langle\rho,\mathcal G\rho\rangle$. Throughout this paper, the free boundary is the
topological boundary of the open superlevel representative, and Lemma \ref{lem:boundary-regular} identifies it with the boundary defined by
positivity of both phase measures in every ball enclosing the point.

The regularity theorem also makes it possible to count distinct energy maximizers and to study their stability. To formulate these additional
conclusions, we define
\begin{equation}\label{eq:maximizer-class}
 \begin{gathered}
 V(m):=\max_{\rho\in\mathcal A_m}\cE(\rho),\qquad
 \M_m:=\{\rho\in\mathcal A_m:\cE(\rho)=V(m)\},\\
 \operatorname*{argmax}_{\omega\in\cA_\lambda}
          \cE_{\mathrm{kin}}(\omega)
 :=\{\lambda\rho:\rho\in\M_{1/\lambda}\}.
 \end{gathered}
\end{equation}
Here we identify patches with their characteristic functions modulo null-measure sets, and write $A\in\M_m$ when $\one_A\in\M_m$. Thus $\#\M_m$ counts
distinct maximizers, whereas the integer $N$ in Theorem \ref{thm:main} counts number of connected-components of energy-maximizing vortex patches. 

The following theorem gives a uniform count over the full mass interval in a smooth convex domain; the boundary may contain line segments.
Smoothness is retained for the endpoint estimates and for the stability assertion. 
\begin{theorem}\label[theorem]{thm:convex-finiteness}
Let $D\subset\R^2$ be a smooth convex bounded domain.
There exists an integer $N_D\geq1$, depending only on $D$, such that
\[
 1\leq\#\M_m\leq N_D\qquad(0<m<|D|).
\]
Moreover, $\#\M_m=1$ for all sufficiently small $m>0$ and all
$m<|D|$ sufficiently close to $|D|$.
\end{theorem}

The proof combines boundary regularity with a nodal argument and a
bilinear Pohozaev identity, followed by analytic reduction and
compactness. It does not require nondegeneracy of the full
mass-constrained second variation. Small-mass uniqueness follows from
\cite{CaoGuoPengYan}, while uniqueness near full
mass is derived from a boundary collar estimate.

\begin{remark}\label[remark]{rem:convex-jordan-finiteness}
For every bounded convex domain $D\subset\R^2$ whose boundary is a rectifiable Jordan curve and every fixed
$m\in(0,|D|)$, there is an integer $N_{D,m}<\infty$ such that
\[
 1\leq\#\M_m\leq N_{D,m}.
\]
No additional regularity assumption on the boundary is required. The proof in Section \ref{end:section} also gives
\[
 \sup_{m\in I}\#\M_m<\infty
 \qquad I \subset \subset(0,|D|) .
\]
The endpoint estimates used for the mass-independent bound in
Theorem \ref{thm:convex-finiteness} require a smooth fluid boundary.
\end{remark}

The finiteness of the set of energy maximizers implies that each maximizer is isolated and hence is a strict local maximizer, with respect to the $L^p$ norm for each $1< p<\infty$, of the kinetic energy on its rearrangement class. Burton's stability theorem \cite[Theorem~5]{Burton2005} then yields the following stability result.
 
 \begin{corollary}\label[corollary]{cor:individual-stability}
Under the assumptions of \cref{thm:convex-finiteness}, let
$\lambda>|D|^{-1}$ and let $\omega_*$ maximize
$\cE_{\mathrm{kin}}$ on $\cA_\lambda$. For every $1< p<\infty$
and $\varepsilon>0$, there is $\delta>0$ such that any Yudovich
solution $\omega(t)$ of \eqref{eq:euler-vorticity} 
satisfies
\[
 \|\omega_0-\omega_*\|_{L^p(D)}<\delta
 \quad\Longrightarrow\quad
 \sup_{t\geq0}\|\omega(t)-\omega_*\|_{L^p(D)}<\varepsilon.
\]
\end{corollary}

The remainder of the paper is organized as follows.
Section \ref{sec:preliminaries} fixes the superlevel representation and the
elliptic estimates. The variational and geometric arguments in
Section \ref{sec:technical} lead to the proof of the main regularity theorem
in Section \ref{S:Proof}. The Dirichlet zero-mode analysis in Section
\ref{ker:section} and the analytic reduction in Section \ref{an:section}
yield fixed-mass finiteness and uniform bounds on compact mass intervals for smooth convex domains. Section \ref{end:section} extends these
counts to arbitrary bounded convex domains, then establishes the endpoint estimates under the smooth boundary assumption. The appendix records conditional extensions to general smooth simply connected domains and identities for smooth families of maximizing patches.

\section{Preliminaries}\label{sec:preliminaries}

Throughout \cref{sec:preliminaries,sec:technical,S:Proof}, $D$ is a bounded
Jordan domain. We first specify the Dirichlet realization and its
boundary continuity, which will keep the maximizing patch away from
$\partial D$ without requiring boundary differentiability. Set
$H^{-1}(D)=(H_0^1(D))^*$, where $H_0^1(D)$ is the closure of
$C_c^\infty(D)$ in $H^1(D)$. Poincar\'e's inequality and the weak
Dirichlet problem define $\mathcal Gf\in H_0^1(D)$ by
\begin{equation}\label{eq:weak-dirichlet}
 \int_D\nabla\mathcal Gf\cdot\nabla z\,\mathrm d\mathbf x
 =\langle f,z\rangle\qquad(z\in H_0^1(D)).
\end{equation}
The Dirichlet Green function is represented as
\[
 G_D(x,y)=\Phi(x-y)-H_D(x,y),\quad
 \Phi(z)=-\frac{1}{2\pi}\log|z|,
\]
where $H_D$ is the regular part of the Green function and
$R_D(x):=H_D(x,x)$ is the Robin function. For $f\in H^{-1}(D)$,
\begin{equation}\label{eq:Hminus-def}
 \Hmone{f}^{2}:=\langle f,\mathcal Gf\rangle
 =\int_D|\nabla\mathcal Gf|^2\,\mathrm d\mathbf x
 =\sup_{0\neq v\in H_0^1(D)}
 \frac{\langle f,v\rangle^2}{\displaystyle\int_D|\nabla v|^2\,\mathrm d\mathbf x}.
\end{equation}
To verify continuity up to the possibly rough fluid boundary, let
$\mathbb D=\{z\in\mathbb C:|z|<1\}$. The Riemann mapping theorem and
Carath\'eodory's theorem \cite{Pommerenke} give a conformal map
\[
 F:\mathbb D\longrightarrow D,
 \qquad F:\overline{\mathbb D}\longrightarrow\overline D
 \quad\text{a homeomorphism}.
\]
For $f\in L^\infty(D)$, extend $f$ by zero and define its Newtonian
potential and pullback by
\[
 N_f(x):=\int_D\Phi(x-y)f(y)\,\mathrm dy,\qquad
 w_f=N_f\circ F.
\]
The potential is $C^{1,\alpha}$ on bounded subsets of $\R^2$ for every
$0<\alpha<1$. Conformal invariance of the Dirichlet integral therefore
gives
\[
 w_f\in H^1(\mathbb D)\cap C(\overline{\mathbb D}),\qquad
 \int_{\mathbb D}|\nabla w_f|^2
 =\int_D|\nabla N_f|^2<\infty.
\]
Let $h_f$ be the harmonic extension of its boundary values. The Poisson
formula and the Dirichlet principle on the disk yield
\[
 \begin{gathered}
 h_f(z)=\frac1{2\pi}\int_0^{2\pi}
       \frac{1-|z|^2}{|e^{it}-z|^2}\,w_f(e^{it})\,\mathrm dt,\qquad
 h_f\in H^1(\mathbb D)\cap C(\overline{\mathbb D}),\\
 w_f-h_f\in H_0^1(\mathbb D),\qquad
 \int_{\mathbb D}|\nabla h_f|^2
 \leq\int_{\mathbb D}|\nabla w_f|^2.
 \end{gathered}
\]
Approximating $w_f-h_f$ by functions in $C_c^\infty(\mathbb D)$,
then composing with $F^{-1}$, preserves convergence of the Dirichlet
integrals; Poincar\'e's inequality on $D$ also gives convergence in
$L^2(D)$. Hence
\[
 (w_f-h_f)\circ F^{-1}\in H_0^1(D),\qquad
 -\Delta\bigl(N_f-h_f\circ F^{-1}\bigr)=f.
\]
Uniqueness in \eqref{eq:weak-dirichlet} identifies the weak solution
with this continuous potential:
\begin{equation}\label{eq:jordan-dirichlet}
 \begin{gathered}
 \mathcal Gf=N_f-h_f\circ F^{-1}
 \in H_0^1(D)\cap C(\overline D)\cap W^{2,p}_{\mathrm{loc}}(D)
 \quad(1<p<\infty),\\
 \mathcal Gf|_{\partial D}=0,\qquad
 |\mathcal Gf|\leq\|f\|_{L^\infty(D)}\tau,
 \qquad \tau:=\mathcal G\one_D\in C(\overline D).
 \end{gathered}
\end{equation}
The last bound follows from the weak comparison principle. Interior
elliptic estimates consequently give
\begin{equation}\label{eq:interior-elliptic}
 \|\mathcal Gf\|_{W^{2,p}(D')}
 \leq C_{D',D,p}\|f\|_{L^\infty(D)}
 \qquad(\overline{D'}\subset D).
\end{equation}
We next verify that the Dirichlet stream function enforces
\eqref{eq:weak-impermeability} on the rough domain. Given
$\psi\in H_0^1(D)$, choose
\[
 \psi_k\in C_c^\infty(D),\qquad
 \psi_k\longrightarrow\psi\quad\text{in }H^1(D).
\]
For every $\varphi\in H^1(D)$, integration by parts is confined to the
compact support of $\psi_k$, and therefore
\[
 \int_D\nabla^\perp\psi_k\cdot\nabla\varphi\,\mathrm d\mathbf x
 =-\int_D\varphi\,\nabla\cdot(\nabla^\perp\psi_k)\,\mathrm d\mathbf x=0.
\]
Passing to the limit gives the required identity without a boundary
trace of the velocity:
\begin{equation}\label{eq:stream-weak-slip}
 \left|\int_D\nabla^\perp\psi\cdot\nabla\varphi\,\mathrm d\mathbf x\right|
 \leq\|\nabla(\psi-\psi_k)\|_{L^2(D)}
      \|\nabla\varphi\|_{L^2(D)}\longrightarrow0.
\end{equation}
In particular, $\mathbf v=\nabla^\perp\mathcal G\omega$ satisfies
\eqref{eq:weak-impermeability} for every bounded vorticity $\omega$.
Testing with $\varphi\in C_c^\infty(D)$ also yields
\[
 \nabla\cdot\mathbf v=0\quad\text{in }\mathcal D'(D).
\]

Let $\Hh^1$ denote the one-dimensional Hausdorff measure and let $\mu\!\llcorner E$ denote restriction to $E$, see for example \cite{SimonGMT}. We note that
 $d\Hh^1=ds$ for regular curves.
On smooth boundaries, $n$ denotes the outward unit normal of $D$
and $\nu$ that of a patch. For a smooth compact curve $\Gamma\subset D$
and $\phi\in L^2(\Gamma)$, the notation
$\phi\,ds|_\Gamma$ denotes the distribution
$z\mapsto\int_\Gamma\phi z\,ds$.
For a measurable $E\subset D$, we define
\begin{equation}\label{eq:area-density}
 \theta_E(\mathbf x,r):=\frac{|E\cap B_r(\mathbf x)|}{\pi r^2},
 \qquad B_r(\mathbf x)\subset D,
\end{equation}
and set
\begin{equation}\label{eq:density-classes}
 E^{(t)}:=\left\{\mathbf x\in D:
 \lim_{r\downarrow0}\theta_E(\mathbf x,r)=t\right\},\qquad t\in\{0,1\}.
\end{equation}
The measure-theoretic boundary is defined as:
\begin{equation}\label{eq:mt-boundary}
 \partial_{\mathrm{mt}}E:=D\setminus(E^{(0)}\cup E^{(1)}).
\end{equation}
Since $\theta_{D\setminus E}=1-\theta_E$, one has
\begin{equation}\label{eq:density-boundary-inclusions}
 \partial_{\mathrm{mt}}E
 =\left\{\limsup_{r\downarrow0}\theta_E>0,\quad
          \limsup_{r\downarrow0}\theta_{D\setminus E}>0\right\}
 \subset\partial E.
\end{equation}
Here $\partial E$ is the topological boundary of the chosen representative. The density classes are unchanged by null-set modifications of $E$.

The following proposition selects an open representative of each
maximizing patch whose boundary lies strictly inside $D$, allowing us to
study the free boundary through interior estimates for the stream
function. It extends the existence and threshold argument of
\cite[Theorem~2.1 and Corollary~2.3]{TurkingtonI} to bounded Jordan
domains; the additional boundary-continuity input is
\eqref{eq:jordan-dirichlet}.
\begin{proposition}\label[proposition]{prop:turkington}
Let $D\subset\R^2$ be a bounded Jordan domain, and
let $\lambda>|D|^{-1}$. The energy $\mathcal E_{\mathrm{kin}}$ attains its
maximum on $\cA_\lambda$. Every maximizer satisfies
\begin{equation}\label{eq:turkington-patch}
 \omega_\lambda=\lambda\one_{A_\lambda}\quad\text{a.e.},\qquad
 A_\lambda=\{\mathcal G\omega_\lambda>\mu_\lambda\}\subset D,\qquad
 |A_\lambda|=\lambda^{-1},
\end{equation}
where $\mu_\lambda>0$ is uniquely determined by $\omega_\lambda$.
Moreover, $\overline{A_\lambda}\subset D$.
Each such maximizer is a weak solution of the steady Euler equation.
\end{proposition}

\begin{proof}
Existence follows from weak-* compactness once continuity of the energy
has been established: the assumption $\lambda>|D|^{-1}$ makes the
constant density admissible, while the mass and box constraints are
preserved under weak-* limits:
\[
 \begin{gathered}
 |D|^{-1}\one_D\in\cA_\lambda,\qquad
 \sup_{\omega\in\cA_\lambda}\|\omega\|_{L^\infty}\leq\lambda,\\
 \omega_n\in\cA_\lambda,\quad\omega_n\overset{*}{\rightharpoonup}\omega
 \ \Longrightarrow\ \omega\in\cA_\lambda.
 \end{gathered}
\]
To pass the quadratic energy to the limit, we use compactness of
$H_0^1(D)\hookrightarrow L^2(D)$, which holds for every bounded domain.
For $\Psi_n:=\mathcal G\omega_n$ and $\Psi:=\mathcal G\omega$, the weak
Dirichlet equation, Poincar\'e's inequality, and uniqueness give
\[
 \begin{gathered}
 \sup_n\|\nabla\Psi_n\|_{L^2(D)}
 \leq C_D\lambda|D|^{1/2},\qquad
 \Psi_n\rightharpoonup\Psi\quad\text{in }H_0^1(D),\\
 \Psi_n\longrightarrow\Psi\quad\text{in }L^2(D),\\
 \|\nabla(\Psi_n-\Psi)\|_{L^2(D)}^2
 =\int_D(\omega_n-\omega)(\Psi_n-\Psi)\,\mathrm d\mathbf x
 \leq2\lambda|D|^{1/2}\|\Psi_n-\Psi\|_{L^2(D)}\longrightarrow0.
 \end{gathered}
\]
Pairing this strong convergence with the bounded densities gives
\[
 \begin{aligned}
 2\bigl(\mathcal E_{\mathrm{kin}}(\omega_n)
       -\mathcal E_{\mathrm{kin}}(\omega)\bigr)
 &=\int_D\omega_n(\Psi_n-\Psi)\,\mathrm d\mathbf x
   +\int_D(\omega_n-\omega)\Psi\,\mathrm d\mathbf x
 \longrightarrow0.
 \end{aligned}
\]
For a maximizing sequence, weak-* sequential compactness yields
\[
 \omega_{n_k}\overset{*}{\rightharpoonup}\omega\in\cA_\lambda,
 \qquad
 \cE_{\mathrm{kin}}(\omega)
 =\lim_{k\to\infty}\cE_{\mathrm{kin}}(\omega_{n_k})
 =\sup_{\cA_\lambda}\cE_{\mathrm{kin}}.
\]

For a fixed maximizer $\omega$, write $\Psi=\mathcal G\omega$ and use
convexity of $\cA_\lambda$ to take one-sided variations towards any
$\eta\in\cA_\lambda$:
\[
 0\geq\left.\frac{\mathrm d}{\mathrm dt}\right|_{t=0+}
 \mathcal E_{\mathrm{kin}}\bigl(\omega+t(\eta-\omega)\bigr)
 =\int_D\Psi(\eta-\omega)\,\mathrm d\mathbf x.
\]
Consequently, $\omega$ also maximizes the linear functional
$\eta\mapsto\int_D\Psi\eta$, so the bathtub principle
\cite[Theorem~1.14]{LiebLoss} supplies a competitor supported on the upper
level sets of $\Psi$: there are $\mu\in\R$ and $\beta\in[0,\lambda]$ such that
\[
 S:=\{\Psi=\mu\},\qquad
 \eta_*:=\lambda\one_{\{\Psi>\mu\}}+\beta\one_S\in\cA_\lambda.
\]
The mass constraint removes the constant $\mu$, while the pointwise
sign of the remaining integrand forces saturation away from the threshold:
\[
 0\geq\int_D\Psi(\eta_*-\omega)\,\mathrm d\mathbf x
 =\int_D(\Psi-\mu)(\eta_*-\omega)\,\mathrm d\mathbf x\geq0,
\]
and hence
\[
 (\lambda-\omega)\one_{\{\Psi>\mu\}}=0,\qquad
 \omega\one_{\{\Psi<\mu\}}=0\quad\text{a.e.}
\]

To exclude vorticity on the possible plateau $S$, we use the Sobolev
level-set identity \cite[Theorem~4.4(iv)]{EvansGariepy}:
\[
 f\in W^{1,p}(D),\ c\in\R
 \quad\Longrightarrow\quad
 \one_{\{f=c\}}\nabla f=0\quad\text{a.e.}
\]
Indeed, for $\chi\in C_c^\infty((-1,1))$, $0\leq\chi\leq1$, $\chi(0)=1$,
\[
 \begin{gathered}
 H_\varepsilon(t):=\int_0^{t-c}\chi(s/\varepsilon)\,\mathrm ds,
 \qquad H_\varepsilon(f)\longrightarrow0\quad\text{in }L^p(D),\\
 \nabla H_\varepsilon(f)
 =\chi((f-c)/\varepsilon)\nabla f
 \longrightarrow\one_{\{f=c\}}\nabla f
 \quad\text{in }L^p(D),\\
 \nabla H_\varepsilon(f)\longrightarrow0\quad\text{in }\mathcal D'(D)
 \quad\Longrightarrow\quad\one_{\{f=c\}}\nabla f=0\quad\text{a.e.}
 \end{gathered}
\]
Since $\Psi\in W^{2,p}_{\mathrm{loc}}(D)$, we apply this identity on
relatively compact subdomains, first to $\Psi$ and then to each
$\partial_i\Psi$. Exhausting $D$ eliminates its second derivatives on
the plateau. With $S_i:=\{\partial_i\Psi=0\}$, $i=1,2$,
\[
 \begin{gathered}
 \one_S\nabla\Psi=0,\qquad |S\setminus S_i|=0,
 \qquad \one_{S_i}\nabla\partial_i\Psi=0\quad\text{a.e.},\\
 \one_S D^2\Psi=0
 \quad\Longrightarrow\quad
 \one_S\omega=-\one_S\Delta\Psi=0\quad\text{a.e.}
 \end{gathered}
\]
Therefore
\[
 \omega=\lambda\one_{\{\Psi>\mu\}}\quad\text{a.e.},
 \qquad |\{\Psi>\mu\}|=\lambda^{-1}.
\]

The threshold must be positive, which in turn separates the patch from
$\partial D$. Indeed, the strong maximum principle gives $\Psi>0$ in
$D$ because $\omega\geq0$ and $\int_D\omega=1$, whereas, with
$M:=\max_{\overline D}\Psi$, either alternative outside $(0,M)$
contradicts the mass constraint.
Therefore
\[
 0<\mu<M,\qquad
 \overline{\{\Psi>\mu\}}
 \subset\{\Psi\geq\mu\}\subset D.
\]
Uniqueness follows from strict decrease of the upper-level-set area
on $(0,M)$, since connectedness of $D$, continuity of $\Psi$, and
$\Psi|_{\partial D}=0$ imply
\[
 \begin{gathered}
 0<a<b<M\quad\Longrightarrow\quad
 \{a<\Psi<b\}\text{ is nonempty and open},\\
 |\{\Psi>a\}|-|\{\Psi>b\}|
 \geq|\{a<\Psi<b\}|>0.
 \end{gathered}
\]
Thus $|\{\Psi>\mu_1\}|=|\{\Psi>\mu_2\}|=\lambda^{-1}$ implies
$\mu_1=\mu_2$.

The superlevel-set representation also gives the steady Euler equation
without requiring regularity of the patch boundary. Indeed, the Sobolev
chain rule yields
\[
 \nabla^\perp(\Psi-\mu)_+
 =\one_{\{\Psi>\mu\}}\nabla^\perp\Psi\quad\text{a.e.},
\]
and therefore, for every $\varphi\in C_c^\infty(D)$,
\[
 \int_D\omega\,\nabla^\perp\Psi\cdot\nabla\varphi\,\mathrm d\mathbf x
 =\lambda\int_D\nabla^\perp(\Psi-\mu)_+\cdot\nabla\varphi
       \,\mathrm d\mathbf x=0,
\]
since the distributional divergence of a perpendicular gradient is zero.
\end{proof}

Set $m=\lambda^{-1}$ and $\rho=\omega/\lambda$. The normalized problem is
\begin{equation}\label{eq:normalized-problem}
 \max_{\rho\in\mathcal A_m}\cE(\rho),\qquad
 \cE(\rho)=\frac12\int_D\rho\mathcal G\rho\,\mathrm d\mathbf x.
\end{equation}
By \cref{prop:turkington}, every energy-maximizing vortex patch
has the superlevel-set representation
\begin{equation}\label{eq:maximizing-patch}
 \rho=\one_A,\qquad \psi=\mathcal G\one_A,\qquad
 u=\psi-q,\qquad A=\{u>0\}\subset D,\qquad
 q=\mu_\lambda/\lambda>0.
\end{equation}
All subsequent boundaries refer to this open representative. The
positive threshold and the zero Dirichlet boundary values of $\psi$
give the compactness that will be used in the boundary arguments:
\[
 \overline A\subset\{\psi\geq q\}\subset D,\qquad
 \operatorname{dist}(\overline A,\partial D)>0.
\]
The normalized stream function satisfies
\begin{equation}\label{eq:PDE-basic}
 -\Delta u=\one_A=\one_{\{u>0\}}\quad\text{in }D,
 \qquad u=-q\quad\text{on }\partial D,
\end{equation}
where the boundary value is attained continuously and
$u+q\in H_0^1(D)$. By \eqref{eq:jordan-dirichlet} and
\eqref{eq:interior-elliptic},
\begin{equation}\label{eq:base-regularity}
 u\in C(\overline D)\cap W^{2,p}_{\mathrm{loc}}(D)
      \cap C^{1,\alpha}_{\mathrm{loc}}(D)
 \qquad(1<p<\infty,\ 0<\alpha<1).
\end{equation}

We shall use the following modulus to control the gradient near a possible singular point on the free boundary, without assuming regularity of $\partial A$.
At $\mathbf x_0\in\partial A$ with $\nabla u(\mathbf x_0)=0$,
\begin{equation}\label{eq:Zygmund}
 |\nabla u(\mathbf x)|\leq C|\mathbf x-\mathbf x_0|
 \left(1+\log\frac{r_0}{|\mathbf x-\mathbf x_0|}\right),
 \qquad0<|\mathbf x-\mathbf x_0|<r_0/2,
\end{equation}
where $B_{r_0}(\mathbf x_0)\subset D$. To isolate the singular part of the
potential, choose
$\chi\in C_c^\infty(B_{r_0}(\mathbf x_0))$, $0\leq\chi\leq1$, $\chi=1$ on
$B_{3r_0/4}(\mathbf x_0)$, and write
\[
 v=\Phi*(\chi\one_A),\qquad u=v+h,
 \qquad\Delta h=0\text{ in }B_{3r_0/4}(\mathbf x_0).
\]
The harmonic remainder has a Lipschitz gradient on the smaller ball.
For the potential, split at twice the separation
$\delta=|\mathbf x-\mathbf x_0|<r_0/2$: use integrability of the gradient
kernel in the near region and the mean-value theorem in the far region.
The bounds
$|\nabla\Phi(\mathbf z)|\leq C|\mathbf z|^{-1}$ and
$|D^2\Phi(\mathbf z)|\leq C|\mathbf z|^{-2}$ give
\[
 \begin{aligned}
 |\nabla v(\mathbf x)-\nabla v(\mathbf x_0)|
 &\leq C\int_{B_{2\delta}(\mathbf x_0)}
 \left(\frac1{|\mathbf x-\mathbf y|}
       +\frac1{|\mathbf x_0-\mathbf y|}\right)\,\mathrm d\mathbf y\\
 &\quad+C\delta\int_{B_{r_0}(\mathbf x_0)\setminus B_{2\delta}(\mathbf x_0)}
                \frac{\mathrm d\mathbf y}{|\mathbf x_0-\mathbf y|^2}\\
 &\leq C\delta+C\delta\int_{2\delta}^{r_0}\frac{\mathrm dt}{t}
 \leq C\delta\left(1+\log\frac{r_0}{\delta}\right),\\
 |\nabla h(\mathbf x)-\nabla h(\mathbf x_0)|
 &\leq\|D^2h\|_{L^\infty(B_{r_0/2}(\mathbf x_0))}\delta\leq C\delta.
 \end{aligned}
\]
Note that $\nabla u(\mathbf x_0)=0$, equation \eqref{eq:Zygmund} is then proved.

The following lemma relates the topological boundary to the two phases without assuming any regularity condition on the boundary, and then establishes local
analyticity wherever the stream-function gradient is nonzero. Together,
these conclusions reduce the remaining regularity problem to critical points of $u$ on $\partial A$.

\begin{lemma}\label[lemma]{lem:boundary-regular}
For an energy-maximizing vortex patch $\one_A$ represented by
\eqref{eq:maximizing-patch}, with supports taken relative to $D$,
\begin{equation}\label{eq:top-fs}
 \partial A
 =\operatorname{supp}_D(\one_A\,\mathrm d\mathbf x)
  \cap\operatorname{supp}_D(\one_{D\setminus A}\,\mathrm d\mathbf x).
\end{equation}
Equivalently, for $\mathbf x\in D$,
\[
 \mathbf x\in\partial A
 \quad\Longleftrightarrow\quad
 \begin{cases}
 |A\cap B_r(\mathbf x)|>0,\\
 |B_r(\mathbf x)\setminus A|>0
 \end{cases}
 \quad\text{for every }0<r<\dist(\mathbf x,\partial D).
\]
If $\mathbf x_0\in\partial A$ and $\nabla u(\mathbf x_0)\neq0$, then
$\partial A$ is a single real-analytic arc near $\mathbf x_0$.
\end{lemma}

\begin{proof}
The boundary identity uses openness of $A=\{u>0\}$ together with the
equation for $u$. If $\mathbf x\in D\setminus\partial A$, a sufficiently
small ball lies in one phase, excluding $\mathbf x$ from the intersection
of the two supports. Conversely, given $\mathbf x\in\partial A$ and
$B_r(\mathbf x)\subset D$, choose
$\mathbf y\in A\cap B_{r/2}(\mathbf x)$; openness then supplies
$\varepsilon\in(0,r/2)$ such that
\[
 B_\varepsilon(\mathbf y)\subset A\cap B_r(\mathbf x)
 \quad\Longrightarrow\quad
 |A\cap B_r(\mathbf x)|\geq\pi\varepsilon^2>0.
\]
No regularity of $\partial A$ is needed for this inclusion.

To prove positivity of the other phase, suppose
$|B_r(\mathbf x)\setminus A|=0$. Continuity of $u$ and
$u>0$ a.e. in the ball imply
\[
 u\geq0\quad\text{on }\overline{B_r(\mathbf x)},\qquad
 -\Delta u=1\quad\text{in }B_r(\mathbf x),\qquad
 u(\mathbf x)=0.
\]
The disk torsion function
\[
 v(\mathbf z):=\frac{r^2-|\mathbf z-\mathbf x|^2}{4},\qquad
 -\Delta v=1\quad\text{in }B_r(\mathbf x),\qquad
 v=0\quad\text{on }\partial B_r(\mathbf x)
\]
then gives a contradiction by the maximum principle:
\[
 \begin{gathered}
 -\Delta(u-v)=0\quad\text{in }B_r(\mathbf x),\qquad
 u-v\geq0\quad\text{on }\partial B_r(\mathbf x),\\
 u\geq v\quad\text{in }B_r(\mathbf x)
 \quad\Longrightarrow\quad
 0=u(\mathbf x)\geq v(\mathbf x)=\frac{r^2}{4}>0.
 \end{gathered}
\]
This establishes \eqref{eq:top-fs}, ensuring that both phases have
positive measure on every fixed ball centered at the boundary; the
uniform half-density estimate will follow later from
\cref{cor:boundary-density}.

At a noncritical point, the implicit function theorem identifies the
boundary with a single $C^{1,\alpha}$ zero-level arc, whose regularity
can be improved to real analyticity by applying
\cite[Lemma~3]{ChanilloGrieserKurata} to
\[
 \Delta u=h(u,\nabla u),\qquad h(s,\mathbf p)=-\one_{\{s>0\}},\qquad
 \partial_{\mathbf p}^{\boldsymbol\alpha}h=0
 \quad(|\boldsymbol\alpha|\geq1).
\]
The hypotheses of that lemma hold since
$\|h\|_{L^\infty}=1$ and $h(s,\cdot)$ is entire; analyticity in the gradient
variable is sufficient, despite the jump in $s$.
\end{proof}

Define
\begin{equation}\label{eq:regular-set}
 \mathcal R:=\partial A\cap\{|\nabla u|>0\}.
\end{equation}
By \cref{lem:boundary-regular}, every $\mathbf x\in\mathcal R$ has a
neighborhood $U\subset D$ and a real-analytic embedding
$\gamma:(-1,1)\to U$ such that
\[
 \inf_U|\nabla u|>0,\qquad
 U\cap\partial A=U\cap\mathcal R=\gamma((-1,1)).
\]
Thus $\mathcal R$ is an embedded real-analytic one-dimensional manifold
without boundary, and $U$ meets exactly one component of $\mathcal R$.

\section{Variational estimates and boundary geometry}\label{sec:technical}

Throughout, $\one_A$ is an energy-maximizing vortex patch, and
$A=\{u>0\}$ is its open vortex core, as chosen in
\cref{sec:preliminaries}.

We begin by quantifying how closely the patch is approximated by a
half-plane near a noncritical boundary point, keeping explicit the
dependence on $|\nabla u(\mathbf p)|$ so that the density estimate can
later be made uniform once critical boundary points have been excluded.

\begin{lemma}\label[lemma]{lem:regular-density}
Let $0<\alpha<1$, $\mathbf p\in\partial A$, and
\[
 g_{\mathbf p}:=|\nabla u(\mathbf p)|>0,\qquad
 \mathbf e_{\mathbf p}:=\nabla u(\mathbf p)/g_{\mathbf p},\qquad
 B_{r_0}(\mathbf p)\subset D,\qquad
 [\nabla u]_{C^{0,\alpha}(B_{r_0}(\mathbf p))}\leq L.
\]
For $0<r<r_0$, set
\[
 A_{\mathbf p,r}:=\{\mathbf y\in B_1:\mathbf p+r\mathbf y\in A\},
 \qquad
 H_{\mathbf p}^+:=\{\mathbf y\in B_1:\mathbf e_{\mathbf p}\cdot\mathbf y>0\}.
\]
Then
\begin{equation}\label{eq:regular-density-estimate}
 |A_{\mathbf p,r}\mathbin{\triangle}H_{\mathbf p}^+|
 \leq\frac{4L}{g_{\mathbf p}}r^\alpha,
 \qquad
 \left|\theta_A(\mathbf p,r)-\frac12\right|
 \leq\frac{4L}{\pi g_{\mathbf p}}r^\alpha.
\end{equation}
\end{lemma}

\begin{proof}
Rescaling to $B_1$ reduces the density estimate to comparison with the
half-disk determined by $\nabla u(\mathbf p)$. Since multiplication by
$r^{-1}>0$ preserves the positivity set, Taylor's formula and
$u(\mathbf p)=0$ give
\[
 \begin{aligned}
 \left|\frac{u(\mathbf p+r\mathbf y)}r
       -\nabla u(\mathbf p)\cdot\mathbf y\right|
 &=\left|\int_0^1
 [\nabla u(\mathbf p+tr\mathbf y)-\nabla u(\mathbf p)]
 \cdot\mathbf y\,\mathrm dt\right|\\
 &\leq\frac{L}{1+\alpha}r^\alpha|\mathbf y|^{1+\alpha}
 \leq Lr^\alpha\qquad(\mathbf y\in B_1).
 \end{aligned}
\]
The signs of $u(\mathbf p+r\mathbf y)$ and
$\mathbf e_{\mathbf p}\cdot\mathbf y$ agree wherever the magnitude of the
linear term exceeds the error. Thus, with $a=Lr^\alpha/g_{\mathbf p}$,
\[
 \begin{gathered}
 A_{\mathbf p,r}\mathbin{\triangle}H_{\mathbf p}^+
 \subset S_a:=\{\mathbf y\in B_1:
                  |\mathbf e_{\mathbf p}\cdot\mathbf y|\leq a\},\\
 |S_a|=\int_{-\min\{a,1\}}^{\min\{a,1\}}
                  2\sqrt{1-s^2}\,\mathrm ds\leq4a,\\
 \left|\theta_A(\mathbf p,r)-\frac12\right|
 =\frac1\pi\bigl||A_{\mathbf p,r}|-|H_{\mathbf p}^+|\bigr|
 \leq\frac1\pi|A_{\mathbf p,r}\mathbin{\triangle}H_{\mathbf p}^+|
 \leq\frac{4a}\pi.
 \end{gathered}
\]
\end{proof}

To turn energy maximization into estimates on the boundary geometry, we
record an exact inequality for competitors of the same mass, which will
then be applied to area-preserving boundary variations.

\begin{proposition}\label[proposition]{prop:VI}
For every $\eta\in L^\infty(D)$ satisfying
$0\leq\eta\leq1$ and $\int_D\eta\,\mathrm d\mathbf x=|A|$,
\begin{equation}\label{eq:VI}
 \int_D u(\eta-\one_A)\,\mathrm d\mathbf x
 +\frac12\Hmone{\eta-\one_A}^2\leq0.
\end{equation}
\end{proposition}

\begin{proof}
Writing $h=\eta-\one_A$, the quadratic structure of the energy and
self-adjointness of $\mathcal G$ give an expansion with no remainder,
while the mass constraint $\int_Dh=0$ eliminates the constant $q$ in
$\mathcal G\one_A=u+q$:
\begin{equation}\label{eq:energy-expansion}
 0\geq\cE(\eta)-\cE(\one_A)
 =\int_D(\mathcal G\one_A)h\,\mathrm d\mathbf x
   +\frac12\langle h,\mathcal Gh\rangle
 =\int_Duh\,\mathrm d\mathbf x+\frac12\Hmone h^2.
\end{equation}
\end{proof}

Applying this inequality to normal variations supported on finitely many
regular boundary pieces yields an estimate that couples all the pieces
through a single total-area constraint.

A \emph{regular arc} is a compact embedded subarc of $\mathcal R$ of positive
length; a \emph{closed regular curve} is a simple closed subcurve of
$\mathcal R$. Let
\[
 \Gamma=\bigsqcup_{\ell=1}^M\Gamma_\ell,\qquad
 g=|\nabla u|\big|_\Gamma,
\]
be a nonempty finite disjoint union of such pieces. The normal
$\boldsymbol\nu$ points out of $A$; thus
\[
 \min_\Gamma g>0,\qquad
 \boldsymbol\nu=-\frac{\nabla u}{|\nabla u|}\Big|_\Gamma,
 \qquad \partial_{\boldsymbol\nu}u=-g.
\]

\begin{proposition}\label[proposition]{prop:shape}
For every $\phi\in L^2(\Gamma)$ with
\begin{equation}\label{eq:zero-flux}
 \int_\Gamma\phi\,\mathrm d\Hh^1=0,
\end{equation}
one has
\begin{equation}\label{eq:shape}
 \left\|\phi\Hh^1\!\llcorner\Gamma\right\|_{\dot H^{-1}_D}^2
 \leq\int_\Gamma g\phi^2\,\mathrm d\Hh^1.
\end{equation}
\end{proposition}

\begin{proof}
Fix $0<\alpha<1$ and initially take $\phi$ smooth, supported away from
arc endpoints, and periodic on closed curves. After slightly extending
each arc within $\mathcal R$, \cref{lem:boundary-regular} provides a
neighborhood of each point of $\Gamma$ in which the boundary is a single
analytic arc. Compactness and pairwise disjointness then give disjoint
tubes, without requiring regularity of the remaining boundary, and in
arclength coordinates we have
\[
 F_\ell(s,r)=\boldsymbol\gamma_\ell(s)+r\boldsymbol\nu_\ell(s),
 \quad A=\{r<0\},\quad
 J_\ell(s,r)=1+\kappa_\ell(s)r\in[1/2,3/2]
 \quad(|r|<\varepsilon),
\]
where $\kappa_\ell$ is the signed curvature determined by the displayed
Jacobian. Although \eqref{eq:zero-flux} cancels the first-order area
change, an additional correction is needed to preserve the area
exactly; for this purpose, choose smooth $\zeta_\ell$ with the same
support conventions and
\[
 B_0:=\sum_\ell\int_{\Gamma_\ell}\zeta_\ell\,\mathrm d\Hh^1\neq0,
 \qquad b_\ell(t,c):=t\phi_\ell+c\zeta_\ell.
\]
The signed area change is
\begin{equation}\label{eq:volume-map}
 \mathcal V(t,c)=\sum_\ell\int_{\Gamma_\ell}
 \left(b_\ell(t,c)+\frac{\kappa_\ell}{2}b_\ell(t,c)^2\right)
 \,\mathrm d\Hh^1.
\end{equation}
The implicit function theorem gives
\[
 \begin{gathered}
 \mathcal V(0,0)=\mathcal V_t(0,0)=0,\qquad
 \mathcal V_c(0,0)=B_0\neq0,\\
 \mathcal V(t,c(t))=0,\qquad c(0)=c'(0)=0,\qquad
 c''(0)=-\frac1{B_0}\sum_\ell
          \int_{\Gamma_\ell}\kappa_\ell\phi_\ell^2\,\mathrm d\Hh^1.
 \end{gathered}
\]
Hence
\begin{equation}\label{eq:c-expansion}
 c(t)=-\frac{t^2}{2B_0}
 \sum_\ell\int_{\Gamma_\ell}\kappa_\ell\phi_\ell^2\,\mathrm d\Hh^1
 +O(t^3).
\end{equation}
With $b_{\ell,t}=b_\ell(t,c(t))$, define $A_t$ by replacing $r<0$ with
$r<b_{\ell,t}(s)$ in the tubes and leaving $A$ unchanged elsewhere.
For $|t|<t_0$, the displacements vanish near the tube ends and satisfy
\[
 \|b_{\ell,t}\|_\infty\leq C_0|t|<\varepsilon/2,
 \qquad |A_t|=|A|.
\]
Thus the correction preserves the first-order normal displacement
$\phi$ while enforcing the area constraint for every small $t$.

Put $h_t=\one_{A_t}-\one_A$ and
$\mu_\phi=\phi\Hh^1\!\llcorner\Gamma$. We identify the quadratic energy
term by proving strong convergence of $t^{-1}h_t$ to $\mu_\phi$ in
$H^{-1}(D)$; weak convergence alone would give only lower semicontinuity
of the squared norm. For $z\in H_0^1(D)$, write
$z_\ell=z\circ F_\ell$. Signed integration gives
\begin{equation}\label{eq:layer-pairing}
 \left\langle\frac{h_t}{t},z\right\rangle
 =\sum_\ell\int_{\Gamma_\ell}\frac1t
 \int_0^{b_{\ell,t}(s)}z_\ell(s,r)(1+\kappa_\ell r)
 \,\mathrm dr\,\mathrm ds.
\end{equation}
The Sobolev slice estimate is
\[
 |z_\ell(s,r)-z_\ell(s,0)|
 \leq |r|^{1/2}
 \left(\int_{-|r|}^{|r|}|\partial_rz_\ell(s,\rho)|^2
 \,\mathrm d\rho\right)^{1/2}.
\]
If $I_\ell$ is the arclength interval, with endpoints identified for a
closed curve, define
\[
 U_\delta:=\bigcup_\ell F_\ell(I_\ell\times(-\delta,\delta)),
 \qquad
 Q_\delta(z)^2:=\sum_\ell\int_{I_\ell}\int_{-\delta}^{\delta}
 |\partial_rz_\ell|^2\,\mathrm dr\,\mathrm ds.
\]
Then
\[
 |U_\delta|\leq3\delta\Hh^1(\Gamma),\qquad
 Q_\delta(z)^2\leq2\|\nabla z\|_{L^2(D)}^2.
\]
Subtracting $\langle\mu_\phi,z\rangle$ in \eqref{eq:layer-pairing} yields
$I_1+I_2+I_3$, where
\begin{align*}
 I_1&=\sum_\ell\int_{\Gamma_\ell}
       \frac{c(t)}t\zeta_\ell z_\ell(s,0)\,\mathrm ds,\\
 I_2&=\sum_\ell\int_{\Gamma_\ell}\frac1t
       \int_0^{b_{\ell,t}}(z_\ell(s,r)-z_\ell(s,0))
       \,\mathrm dr\,\mathrm ds,\\
 I_3&=\sum_\ell\int_{\Gamma_\ell}\frac1t
       \int_0^{b_{\ell,t}}\kappa_\ell r z_\ell(s,r)
       \,\mathrm dr\,\mathrm ds.
\end{align*}
The trace, slice, and Poincar\'e inequalities imply
\begin{align}
 |I_1|&\leq C\frac{|c(t)|}{|t|}\|\nabla z\|_{L^2(D)},\notag\\*
 |I_2|&\leq C|t|^{1/2}Q_{C_0|t|}(z)
       \leq C|t|^{1/2}\|\nabla z\|_{L^2(D)},\notag\\*
 |I_3|&\leq C\|z\|_{L^1(U_{C_0|t|})}
       \leq C|t|^{1/2}\|\nabla z\|_{L^2(D)}.
 \label{eq:layer-errors}
\end{align}
Here $C$ is independent of $t,z$. Consequently,
\begin{equation}\label{eq:strong-Hminus}
 \Hmone{t^{-1}h_t-\mu_\phi}
 \leq C\left(|t|^{1/2}+\frac{|c(t)|}{|t|}\right)
 \leq C|t|^{1/2}\longrightarrow0,
\end{equation}
and therefore
\begin{equation}\label{eq:norm-expansion}
 \Hmone{h_t}^2=t^2\Hmone{\mu_\phi}^2+o(t^2).
\end{equation}

The linear energy term depends on the first normal derivative of $u$.
Since $u=0$ on $\Gamma$ and $\partial_{\boldsymbol\nu}u=-g$,
$u\circ F_\ell(s,r)=-g(s)r+O(|r|^{1+\alpha})$ uniformly on the fixed
family of tubes. Signed integration, valid for both positive and
negative displacements, gives
\begin{align}
 \int_Duh_t\,\mathrm d\mathbf x
 &=\sum_\ell\int_{\Gamma_\ell}\int_0^{b_{\ell,t}}
 (u\circ F_\ell)(s,r)(1+\kappa_\ell r)\,\mathrm dr\,\mathrm ds\notag\\
 &=-\frac12\sum_\ell\int_{\Gamma_\ell}g b_{\ell,t}^2\,\mathrm d\Hh^1
   +O(|t|^{2+\alpha}+|t|^3)\notag\\
 &=-\frac{t^2}{2}\int_\Gamma g\phi^2\,\mathrm d\Hh^1+o(t^2).
 \label{eq:u-expansion}
\end{align}
Since the construction preserves the area exactly, we may apply
\eqref{eq:VI} with $\eta=\one_{A_t}$ and combine the two energy
expansions to obtain
\[
 0\geq\frac{2}{t^2}
       \bigl(\cE(\one_{A_t})-\cE(\one_A)\bigr)
  =\frac{2}{t^2}\int_Duh_t\,\mathrm d\mathbf x
    +\Hmone{t^{-1}h_t}^2 \longrightarrow
 -\int_\Gamma g\phi^2\,\mathrm d\Hh^1+\Hmone{\mu_\phi}^2.
\]

To extend the inequality to $L^2$ displacements, approximate while
retaining the total-flux constraint. Choose smooth
$\psi_n\to\phi$ in $L^2(\Gamma)$ with the preceding support conventions,
and fix such a smooth $\zeta$ with $\int_\Gamma\zeta=1$.
Then
\[
 \phi_n:=\psi_n-\left(\int_\Gamma\psi_n\,\mathrm d\Hh^1\right)\zeta,
 \qquad \int_\Gamma\phi_n=0,\qquad
 \phi_n\to\phi\quad\text{in }L^2(\Gamma).
\]
The trace theorem gives, for $f\in L^2(\Gamma)$,
\[
 \begin{aligned}
 \Hmone{f\Hh^1\!\llcorner\Gamma}
 &=\sup_{\substack{z\in H_0^1(D)\\\|\nabla z\|_{L^2(D)}=1}}
          \left|\int_\Gamma fz\,\mathrm d\Hh^1\right|\\
 &\leq\|f\|_{L^2(\Gamma)}
       \sup_{\substack{z\in H_0^1(D)\\\|\nabla z\|_{L^2(D)}=1}}
       \|z\|_{L^2(\Gamma)}
 \leq C_\Gamma\|f\|_{L^2(\Gamma)}.
 \end{aligned}
\]
Thus
\[
 \Hmone{(\phi_n-\phi)\Hh^1\!\llcorner\Gamma}
 \leq C_\Gamma\|\phi_n-\phi\|_{L^2(\Gamma)}\to0,
 \qquad
 \int_\Gamma g\phi_n^2\,\mathrm d\Hh^1
 \to\int_\Gamma g\phi^2\,\mathrm d\Hh^1.
\]
Taking $t\to0$ for each smooth $\phi_n$ and then $n\to\infty$ using
these two convergences proves \eqref{eq:shape}. Because the constraint
acts on the total flux $\int_\Gamma\phi$, it allows flux to be transferred
between distinct pieces, while the quadratic term retains their Green
interactions:
\[
 \Hmone{\mu_\phi}^2
 =\sum_{\ell,k=1}^M\int_{\Gamma_\ell}\int_{\Gamma_k}
 G_D(x,y)\phi_\ell(x)\phi_k(y)
 \,\mathrm d\Hh^1(y)\,\mathrm d\Hh^1(x).
\]
\end{proof}

The variation estimate can now be expressed as a weighted trace
inequality for test functions in $H_0^1(D)$, in a form suited to
controlling the weighted boundary measure $g^{-1}\,d\Hh^1$ on families of regular
arcs that approach a possible critical point.

\begin{corollary}\label[corollary]{cor:trace}
For $v\in H_0^1(D)$, set
\[
 c_{\Gamma,v}:=
 \frac{\displaystyle\int_\Gamma v/g\,\mathrm d\Hh^1}
      {\displaystyle\int_\Gamma 1/g\,\mathrm d\Hh^1}.
\]
Then
\begin{equation}\label{eq:weighted-trace}
 \int_\Gamma\frac{(v-c_{\Gamma,v})^2}{g}\,\mathrm d\Hh^1
 \leq\int_D|\nabla v|^2\,\mathrm d\mathbf x.
\end{equation}
\end{corollary}

\begin{proof}
Subtracting the common weighted mean makes the displacement below
satisfy the single area constraint. Since $\min_\Gamma g>0$,
the trace of $v$ defines
\[
 c:=c_{\Gamma,v},\qquad
 \phi:=\frac{v-c}g\in L^2(\Gamma),\qquad
 \int_\Gamma\phi=0,\qquad
 X:=\int_\Gamma g\phi^2\,\mathrm d\Hh^1.
\]
Pairing this displacement with $v$ recovers the weighted boundary
variance $X$, to which Dirichlet duality and \cref{prop:shape} apply:
\[
 X=\int_\Gamma(v-c)\phi\,\mathrm d\Hh^1
  =\int_\Gamma v\phi\,\mathrm d\Hh^1
 \leq\Hmone{\phi\Hh^1\!\llcorner\Gamma}\|\nabla v\|_{L^2(D)}
 \leq X^{1/2}\|\nabla v\|_{L^2(D)}.
\]
For $X=0$ the claim is immediate; for $X>0$, division by $X^{1/2}$ gives
$X\leq\|\nabla v\|_{L^2(D)}^2$.
\end{proof}

The weighted trace inequality controls regular boundary pieces
approaching a common point through the following summability criterion,
which will be used to exclude accumulation of regular curves at a
critical point.

\begin{lemma}\label[lemma]{lem:annular}
Let $\mathbf{x}_0\in D$ and
\[
 R_1>R_2>\cdots\downarrow0,
 \qquad B_{R_1}(\mathbf{x}_0)\subset D.
\]
Suppose that $\gamma_j$ are pairwise disjoint compact regular arcs or closed
regular curves with
\[
 \gamma_j\subset
 \overline{B_{R_j}(\mathbf{x}_0)}\setminus B_{R_{j+1}}(\mathbf{x}_0).
\]
Then
\begin{equation}\label{eq:annular-sum}
 \sum_{j=1}^\infty\int_{\gamma_j}\frac{1}{|\nabla u|}
 \,\mathrm d\Hh^1<\infty.
\end{equation}
If $\nabla u(\mathbf{x}_0)=\mathbf0$, there are no $R,c>0$ and regular arcs
satisfying, for every sufficiently large $k$,
\begin{equation}\label{eq:long-dyadic-arcs}
 \gamma_k\subset\overline{B_{2^{-k}R}(\mathbf{x}_0)}
 \setminus B_{2^{-k-1}R}(\mathbf{x}_0),
 \qquad\Hh^1(\gamma_k)\geq c\,2^{-k}R.
\end{equation}
\end{lemma}

\begin{proof}
Set
\[
 g=|\nabla u|,\qquad
 a_j=\int_{\gamma_j}g^{-1}\,\mathrm d\Hh^1\in(0,\infty),
 \qquad A_K=\sum_{j=1}^K a_j.
\]
We apply the weighted trace inequality to two separated families with
equal $g^{-1}\,d\Hh^1$ mass, choosing a cutoff that equals $0$ on one
family and $2$ on the other. Its weighted mean is then $1$, whereas
its Dirichlet energy is independent of the spatial scale in dimension
two, giving the contradiction below if the weighted masses diverge.
Suppose $A_K\to\infty$. Choose
\[
 \begin{gathered}
 V\in C_c^\infty(B_1),\qquad 0\leq V\leq2,\qquad
 V|_{B_{1/2}}=2,\\
 v_K(\mathbf{x})=V\left(\frac{\mathbf{x}-\mathbf{x}_0}{R_{K+1}}\right).
 \end{gathered}
\]
For $M>K+1$ such that $R_M<R_{K+1}/2$,
\[
 v_K|_{\gamma_j}=0\quad(j\leq K),\qquad
 v_K|_{\gamma_j}=2\quad(j\geq M),\qquad
 \sum_{j=M}^\infty a_j=\infty.
\]
Since the tail has infinite weighted mass, we can truncate it to match
the mass $A_K$ of the first $K$ pieces. Let
\[
 L=\min\left\{\ell\geq M:\sum_{j=M}^{\ell}a_j\geq A_K\right\},
 \qquad
 \beta_K=A_K-\sum_{j=M}^{L-1}a_j\in(0,a_L].
\]
For an arclength parametrization
$\boldsymbol\gamma_L:[0,\Hh^1(\gamma_L)]\to\gamma_L$, choose $s_K$ by
\[
 \int_0^{s_K}\frac{\mathrm ds}{g(\boldsymbol\gamma_L(s))}=\beta_K.
\]
The integral is continuous and strictly increasing from $0$ to $a_L$,
so this choice realizes the remaining mass $\beta_K$ exactly. Set
\[
 \begin{gathered}
 T_K=\bigcup_{j=M}^{L-1}\gamma_j\,
       \cup\boldsymbol\gamma_L([0,s_K]),\qquad
 \Gamma_K=\bigcup_{j=1}^K\gamma_j\,\cup T_K,\\
 \int_{T_K}g^{-1}\,\mathrm d\Hh^1=A_K.
 \end{gathered}
\]
Then
\[
 \int_{\Gamma_K}g^{-1}\,\mathrm d\Hh^1=2A_K,\qquad
 \int_{\Gamma_K}v_Kg^{-1}\,\mathrm d\Hh^1=2A_K,\qquad
 c_{\Gamma_K,v_K}=1.
\]
The last piece is a compact subarc or the whole closed curve.
By \cref{cor:trace} and two-dimensional Dirichlet scaling,
\begin{equation}\label{eq:two-AK}
 2A_K
 =\int_{\Gamma_K}\frac{(v_K-1)^2}{g}\,\mathrm d\Hh^1
 \leq\int_D|\nabla v_K|^2\,\mathrm d\mathbf{x}
 =\int_{\R^2}|\nabla V|^2\,\mathrm d\mathbf{x},
\end{equation}
contrary to $A_K\to\infty$. The last equality follows because the
change of variables $\mathbf x=\mathbf x_0+R_{K+1}\mathbf y$ cancels
the factor $R_{K+1}^{-2}$ from $|\nabla v_K|^2$ against the Jacobian
$R_{K+1}^{2}$, so no comparison between consecutive radii is required.

For the second assertion, suppose \eqref{eq:long-dyadic-arcs} holds and set
$\rho_j=2^{-2j}R$. Retaining the even indices separates the arcs by a
positive radial gap, including when consecutive dyadic arcs meet on a
common circle. Then
\[
 \begin{gathered}
 \mathbf x_0\in\overline{\bigcup_j\gamma_{2j}}\subset\partial A,
 \\ 
 \gamma_{2j}\subset
 \overline{B_{\rho_j}(\mathbf{x}_0)}\setminus B_{\rho_{j+1}}(\mathbf{x}_0),
 \qquad\Hh^1(\gamma_{2j})\geq c\rho_j.
 \end{gathered}
\]
The gradient estimate at a critical point converts the lower bound on
arc length into a divergent lower bound on weighted mass. By
\eqref{eq:Zygmund}, for
$j\geq j_0$ with $j_0$ sufficiently large,
\[
 \sup_{\gamma_{2j}}g\leq C\rho_j(1+\log(2r_0/\rho_j))
 \leq C'(j+1)\rho_j,
 \qquad
 \int_{\gamma_{2j}}g^{-1}\,\mathrm d\Hh^1\geq\frac{c'}{j+1}.
\]
Thus
\[
 \infty=\sum_{j\geq j_0}\frac{c'}{j+1}
 \leq\sum_{j\geq j_0}\int_{\gamma_{2j}}g^{-1}\,\mathrm d\Hh^1
 <\infty,
\]
contradicting \eqref{eq:annular-sum}. This argument follows
the device of \cite[Lemmas~3.1--3.2]{ChanilloKenigTo}.
\end{proof}

To recover the entire boundary from its regular part, we first show that
every boundary point can be approximated by regular points.

\begin{lemma}\label[lemma]{lem:density}
The regular set $\mathcal R$ is dense in $\partial A$.
\end{lemma}

\begin{proof}
An interior disk tangent to the patch boundary provides a point with
nonzero gradient, because comparison with its explicit torsion function
gives a derivative bound without requiring regularity at the contact
point. To construct such a point near any given
$\mathbf p\in\partial A$, fix $B_\varepsilon(\mathbf p)\subset D$ and
choose
\[
 \begin{gathered}
 \mathbf a\in A\cap B_{\varepsilon/4}(\mathbf p),\qquad
 \mathbf z\in\partial A,\\
 0<d:=\dist(\mathbf a,\partial A)=|\mathbf z-\mathbf a|
 \leq|\mathbf a-\mathbf p|<\varepsilon/4.
 \end{gathered}
\]
Then
\[
 B_d(\mathbf a)\subset A,\qquad
 \overline{B_d(\mathbf a)}\subset B_{\varepsilon/2}(\mathbf p),\qquad
 v(\mathbf x):=\frac{d^2-|\mathbf x-\mathbf a|^2}{4}.
\]
The maximum principle gives
\[
 \left\{\begin{aligned}
  \Delta(u-v)&=0&&\text{in }B_d(\mathbf a),\\
  u-v&\geq0&&\text{on }\partial B_d(\mathbf a),
 \end{aligned}\right.
 \quad\Longrightarrow\quad u\geq v\quad\text{in }B_d(\mathbf a).
\]
For
$\mathbf n=(\mathbf z-\mathbf a)/d$,
\[
 u(\mathbf z-t\mathbf n)\geq\frac{dt}{2}-\frac{t^2}{4}
 \quad(0<t<d),\qquad
 -\nabla u(\mathbf z)\cdot\mathbf n
 =\lim_{t\downarrow0}\frac{u(\mathbf z-t\mathbf n)}{t}
 \geq\frac d2>0.
\]
Therefore
\[
 \mathcal R\cap B_{\varepsilon/2}(\mathbf p)\ni\mathbf z,
 \qquad \mathbf p\in\overline{\mathcal R}.
\]
Here $\varepsilon>0$ can be arbitrarily small.
\end{proof}

Having established density of the regular set, we use the annular
criterion to show that each of its components closes up inside $D$,
rather than accumulating at a critical boundary point.

\begin{lemma}\label[lemma]{lem:continuation}
Every connected component of $\mathcal R$ is a simple closed real-analytic
curve compactly contained in $D$.
\end{lemma}

\begin{proof}
For a connected component $C$ of $\mathcal R$, accumulation at a critical
point would force $C$ to cross every sufficiently small dyadic annulus,
contrary to \cref{lem:annular}. We make this argument precise using
\[
 \overline C\subset\partial A\subset D.
\]
Suppose $\mathbf p\in\overline C\setminus\mathcal R$, so
$\nabla u(\mathbf p)=\mathbf0$. Choose
\[
 \mathbf z\in C,\qquad B_{\rho_0}(\mathbf p)\subset D,\qquad
 \rho_k=2^{-k}\rho_0.
\]
For all sufficiently large $k$,
\[
 \mathbf y_k\in C\cap B_{\rho_k/2}(\mathbf p),\qquad
 |\mathbf z-\mathbf p|>\rho_k.
\]
Arcwise connectedness of $C$ gives a smooth embedded arc
\[
 \boldsymbol\sigma_k:[0,1]\to C,\qquad
 \boldsymbol\sigma_k(0)=\mathbf y_k,\qquad
 \boldsymbol\sigma_k(1)=\mathbf z.
\]
Set
\[
 \begin{split}
 b_k&=\min\{t:|\boldsymbol\sigma_k(t)-\mathbf p|=\rho_k\},\\
 a_k&=\max\{t\leq b_k:|\boldsymbol\sigma_k(t)-\mathbf p|=\rho_k/2\},\\
 \gamma_k&=\boldsymbol\sigma_k([a_k,b_k]).
 \end{split}
\]
Then $a_k<b_k$ and
\[
 \gamma_k\subset\overline{B_{\rho_k}(\mathbf p)}
                  \setminus B_{\rho_k/2}(\mathbf p),\qquad
 \Hh^1(\gamma_k)
 \geq|\boldsymbol\sigma_k(b_k)-\boldsymbol\sigma_k(a_k)|
 \geq\rho_k/2.
\]
Thus \cref{lem:annular} gives
\[
 \overline C\setminus\mathcal R=\varnothing.
\]

It remains to exclude accumulation on a different regular component.
For $\mathbf p\in\overline C$, choose a coordinate neighborhood $U$ such
that $U\cap\mathcal R$ is one connected arc. Then
\[
 (U\cap\mathcal R)\cap C\neq\varnothing
 \quad\Longrightarrow\quad U\cap\mathcal R\subset C
 \quad\Longrightarrow\quad\mathbf p\in C.
\]
Consequently,
\[
 C=\overline C\subset D,\qquad
 0<\Hh^1(C)<\infty,\qquad \min_C|\nabla u|>0.
\]
The classification of compact connected one-dimensional manifolds
without boundary shows that $C$ is diffeomorphic to $\mathbb S^1$;
together with its embedding and \cref{lem:boundary-regular}, this
identifies $C$ as a real-analytic Jordan curve.
\end{proof}

To control the number of these closed curves, we identify the phase
enclosed by each curve and obtain a lower bound on its weighted boundary
integral that is independent of its size.

\begin{lemma}\label[lemma]{lem:torsion}
For every connected component $\Gamma$ of $\mathcal R$, its Jordan interior
$\Omega_\Gamma=\Jint(\Gamma)$ satisfies
\begin{equation}\label{eq:Jordan-positive}
 \overline{\Omega_\Gamma}\subset D,\qquad \Omega_\Gamma\subset A,
\end{equation}
and
\begin{equation}\label{eq:torsion-bound}
 \int_\Gamma|\nabla u|\,\mathrm d\Hh^1=|\Omega_\Gamma|,
 \qquad\int_\Gamma\frac1{|\nabla u|}\,\mathrm d\Hh^1\geq4\pi.
\end{equation}
\end{lemma}

\begin{proof}
Before testing the equation on the Jordan interior, we use simple
connectedness of $D$ to ensure that this interior lies within the domain
of the equation. Denoting the winding number by
$\operatorname{Ind}(\Gamma,\mathbf y)$, a point
$\mathbf y\in\Omega_\Gamma\setminus D$ would give, by Jordan separation
and contraction of $\Gamma$ in $D$, respectively,
\[
 |\operatorname{Ind}(\Gamma,\mathbf y)|=1,\qquad
 \operatorname{Ind}(\Gamma,\mathbf y)=0,
\]
a contradiction. Therefore
\[
 \Omega_\Gamma\subset D,\qquad
 \overline{\Omega_\Gamma}=\Omega_\Gamma\cup\Gamma\subset D.
\]

Since $u=0$ on $\Gamma$, its negative part
$u^-:=\max\{-u,0\}$ belongs to $H_0^1(\Omega_\Gamma)$ and is therefore
an admissible test function; using $\one_Au^-=0$, which follows from
$A=\{u>0\}$, we obtain
\[
 -\int_{\Omega_\Gamma}|\nabla u^-|^2\,\mathrm d\mathbf x
 =\int_{\Omega_\Gamma}\one_Au^-\,\mathrm d\mathbf x=0,
 \qquad u^-\equiv0.
\]
Since
\[
 u\equiv0\text{ in }\Omega_\Gamma
 \quad\Longrightarrow\quad\nabla u|_\Gamma=\mathbf0
 \quad\Longrightarrow\quad\Gamma\cap\mathcal R=\varnothing,
\]
$u\not\equiv0$, so the strong minimum principle for the nonnegative
superharmonic function $u$ excludes every interior zero and gives
\[
 u>0,\quad-\Delta u=1\quad\text{in }\Omega_\Gamma,
 \qquad\partial_{\boldsymbol\nu}u=-|\nabla u|\quad\text{on }\Gamma,
\]
where $\boldsymbol\nu$ is the outward unit normal. Thus
\[
 |\Omega_\Gamma|
 =-\int_{\Omega_\Gamma}\Delta u\,\mathrm d\mathbf x
 =-\int_\Gamma\partial_{\boldsymbol\nu}u\,\mathrm d\Hh^1
 =\int_\Gamma|\nabla u|\,\mathrm d\Hh^1.
\]
Combining this relation between boundary flux and enclosed area with
the isoperimetric and Cauchy--Schwarz inequalities gives
\[
 4\pi|\Omega_\Gamma|
 \leq\Hh^1(\Gamma)^2
 \leq\left(\int_\Gamma|\nabla u|\,\mathrm d\Hh^1\right)
       \left(\int_\Gamma\frac1{|\nabla u|}\,\mathrm d\Hh^1\right)
 =|\Omega_\Gamma|\int_\Gamma\frac1{|\nabla u|}\,\mathrm d\Hh^1.
\]
Division by $|\Omega_\Gamma|>0$ proves \eqref{eq:torsion-bound}, whose
lower bound $4\pi$ is independent of both the size and shape of the
component.
\end{proof}

The uniform lower bound in \cref{lem:torsion}, together with annular
summability, now rules out infinitely many regular components.

\begin{proposition}\label[proposition]{prop:finiteness}
The regular set $\mathcal R$ has only finitely many connected components.
\end{proposition}

\begin{proof}
If infinitely many components accumulate, either all but finitely many
of these curves leave some fixed neighborhood or every neighborhood
contains infinitely many entire curves. The first possibility produces long
annular arcs; the second contradicts summability through the uniform
torsion bound.
Suppose $\{\Gamma_i\}_{i\geq1}$ are distinct components. Compactness of
$\partial A$ gives, after passage to a subsequence,
\[
 \mathbf x_i\in\Gamma_i,\qquad
 \mathbf x_i\longrightarrow\mathbf x_\infty\in\partial A.
\]
A coordinate neighborhood of a regular point meets only one component;
therefore
\[
 \nabla u(\mathbf x_\infty)=\mathbf0,\qquad
 \dist(\mathbf x_\infty,\Gamma_i)>0\quad\text{for every }i.
\]
We first claim
\begin{equation}\label{eq:entire-curves-claim}
 \#\{i:\Gamma_i\subset B_R(\mathbf x_\infty)\}=\infty
 \qquad\text{whenever }B_R(\mathbf x_\infty)\subset D.
\end{equation}
If \eqref{eq:entire-curves-claim} fails, discarding finitely many curves gives
\[
 \Gamma_i\setminus B_R(\mathbf x_\infty)\neq\varnothing
 \quad(i\geq1),\qquad \mathbf x_i\to\mathbf x_\infty.
\]
Choose distinct $i(k)$ such that
\[
 \Gamma_{i(k)}\cap B_{2^{-k-1}R}(\mathbf x_\infty)\neq\varnothing,
 \qquad\Gamma_{i(k)}\setminus B_R(\mathbf x_\infty)\neq\varnothing.
\]
The hitting-time construction in \cref{lem:continuation} supplies compact
regular arcs
\[
 \gamma_k\subset\Gamma_{i(k)}\cap
 \left(\overline{B_{2^{-k}R}(\mathbf x_\infty)}
       \setminus B_{2^{-k-1}R}(\mathbf x_\infty)\right),
 \qquad\Hh^1(\gamma_k)\geq2^{-k-1}R.
\]
Lemma~\ref{lem:annular} gives a contradiction, proving
\eqref{eq:entire-curves-claim}.

The preceding claim allows us to place entire curves in separated
annuli, with each new curve lying closer to $\mathbf x_\infty$ than
the previously selected one. Choose
$B_{R_1}(\mathbf x_\infty)\subset D$ and inductively select
\[
 \widehat\Gamma_j\subset B_{R_j}(\mathbf x_\infty),\qquad
 d_j=\dist(\mathbf x_\infty,\widehat\Gamma_j),\qquad
 R_{j+1}=d_j/2.
\]
Then
\[
 0<R_{j+1}<R_j/2,\qquad
 \widehat\Gamma_j\subset B_{R_j}(\mathbf x_\infty)
     \setminus\overline{B_{R_{j+1}}(\mathbf x_\infty)}.
\]
The annuli are disjoint, whereas every selected curve contributes at
least $4\pi$ to the weighted integral. By
Lemmas \ref{lem:annular} and \ref{lem:torsion}, we have
\[
\sum_{j=1}^\infty4\pi
 \leq\sum_{j=1}^\infty\int_{\widehat\Gamma_j}
             \frac1{|\nabla u|}\,\mathrm d\Hh^1
 <\infty,
\]
which is a contradiction.
\end{proof}

\section{Proof of the boundary regularity}\label{S:Proof}

Combining density of regular points with finiteness of their closed
components identifies the entire boundary, from which we recover the
positive phase without imposing regularity on the remaining boundary
points.

\begin{proposition}\label[proposition]{prop:exhaustion}
There are finitely many pairwise disjoint real-analytic Jordan curves
$\Gamma_1,\ldots,\Gamma_N$ such that
\begin{equation}\label{eq:boundary-exhaustion}
 \partial A=\bigsqcup_{j=1}^N\Gamma_j,
 \qquad
 A=\bigsqcup_{j=1}^N\Jint(\Gamma_j).
\end{equation}
Every component of $A$ is simply connected, and
\begin{equation}\label{eq:no-critical-boundary}
 \partial A=\mathcal R,\qquad
 \partial A\cap\{\nabla u=\mathbf0\}=\varnothing.
\end{equation}
Moreover,
\begin{equation}\label{eq:full-zero-set}
 \{u=0\}=\partial A,\qquad u<0\quad\text{in }D\setminus\overline A.
\end{equation}
\end{proposition}

\begin{proof}
Since the regular set consists of finitely many compact components,
their union is closed; density therefore leaves no additional boundary
points. More precisely, \cref{lem:continuation,prop:finiteness,lem:density}
give
\[
 1\leq N<\infty,\qquad
 \mathcal R=\bigsqcup_{j=1}^N\Gamma_j\subset D,
 \qquad \overline{\mathcal R}=\mathcal R,\qquad
 \partial A=\overline{\mathcal R}=\mathcal R.
\]
Continuity of the gradient on this compact set then gives
\[
 \partial A\cap\{\nabla u=0\}=\varnothing,\qquad
 \min_{\partial A}|\nabla u|>0.
\]

To recover the patch from these curves, set
$\Omega_j=\Jint(\Gamma_j)$ and recall from \cref{lem:torsion} that
each interior lies in the positive phase. Since two disjoint Jordan
curves can have overlapping interiors only when one lies inside the
other, Jordan separation gives
\[
 \begin{gathered}
 \Omega_j\subset A,\qquad u|_{\Gamma_j}=0,\qquad u|_{\Omega_j}>0,\\
 i\ne j,\quad \Omega_i\cap\Omega_j\ne\varnothing
 \ \Longrightarrow\quad
 \Gamma_i\subset\Omega_j\ \text{or}\ \Gamma_j\subset\Omega_i
 \ \Longrightarrow\quad (\Gamma_i\cup\Gamma_j)\cap A\ne\varnothing.
 \end{gathered}
\]
The identity $(\Gamma_i\cup\Gamma_j)\cap A=\varnothing$ excludes such
nesting, so the interiors are pairwise disjoint; moreover, each
$\Omega_j$ is connected and both open and closed relative to $A$, whence
\[
 \overline{\Omega_j}\cap A=\Omega_j,\qquad
 \Omega_j\text{ is open and connected}
 \quad\Longrightarrow\quad
 \Omega_j\text{ is a connected component of }A.
\]

It remains to exclude a component of $A$ not represented by these
interiors. For any component $A_0$ of $A$, openness and boundedness give
\[
 \varnothing\ne\partial A_0\subset\partial A
 =\bigcup_{j=1}^N\Gamma_j.
\]
Choose $\mathbf x\in\partial A_0\cap\Gamma_j$ and a regular coordinate
neighborhood $U$ of $\mathbf x$, where the positive side agrees with
$\Omega_j$. Consequently,
\[
 \varnothing\ne U\cap A_0\subset U\cap A=U\cap\Omega_j
 \quad\Longrightarrow\quad A_0=\Omega_j,
 \qquad A=\bigsqcup_{j=1}^N\Omega_j.
\]
Each $\Omega_j$ is simply connected by the Jordan--Schoenflies theorem.

To identify the zero set outside the patch, observe that the pairwise
disjoint compact sets $\overline\Omega_j$ lie in $D$ and have smooth
boundaries admitting disjoint exterior collars. A path in $D$ can
therefore be perturbed to meet these boundaries transversely and detoured
through the collars whenever it crosses an interior, proving that
$E:=D\setminus\overline A$ is connected. On this exterior phase,
\[
 \begin{gathered}
 \partial E=\partial D\cup\bigcup_{j=1}^N\Gamma_j,\qquad
 \Delta u=0\ \text{in }E,\qquad
 u|_{\partial D}=-q<0,\qquad u|_{\Gamma_j}=0,\\
 \text{strong maximum principle}
 \quad\Longrightarrow\quad u<0\ \text{in }E
 \quad\Longrightarrow\quad \{u=0\}=\partial A.
 \end{gathered}
\]
\end{proof}

Since the gradient has a positive minimum on the compact boundary of
the fixed patch, the local half-plane approximation in
\cref{lem:regular-density} now yields a uniform density estimate.

\begin{corollary}\label[corollary]{cor:boundary-density}
For each $0<\alpha<1$ there are $C,r_*>0$, depending on the fixed patch,
such that
\begin{equation}\label{eq:uniform-half-density}
 \sup_{\mathbf p\in\partial A}
 \left|\theta_A(\mathbf p,r)-\frac12\right|
 \leq Cr^\alpha\qquad(0<r<r_*).
\end{equation}
In particular,
\begin{equation}\label{eq:half-density}
 \lim_{r\downarrow0}\theta_A(\mathbf p,r)=\frac12
 \qquad(\mathbf p\in\partial A),
\end{equation}
and
\begin{equation}\label{eq:no-endpoint-density}
 \partial A\cap(A^{(0)}\cup A^{(1)})=\varnothing,\qquad
 \partial A=\partial_{\mathrm{mt}}A.
\end{equation}
\end{corollary}

\begin{proof}
Choose $r_0>0$ so small that the compact neighborhood
\[
 K:=\{x\in\R^2:\dist(x,\partial A)\leq2r_0\}\subset D.
\]
By \cref{prop:exhaustion} and the interior regularity in
\eqref{eq:base-regularity},
\[
 g_0:=\min_{\partial A}|\nabla u|>0,\qquad
 L:=[\nabla u]_{C^{0,\alpha}(K)}<\infty.
\]
Choose $0<r_*\leq r_0$. \Cref{lem:regular-density} gives
\[
 \sup_{\mathbf p\in\partial A}
 \left|\theta_A(\mathbf p,r)-\frac12\right|
 \leq\frac{4L}{\pi g_0}r^\alpha\qquad(0<r<r_*).
\]
The limiting density $1/2$ places every topological boundary point in
the measure-theoretic boundary, while
\eqref{eq:density-boundary-inclusions} supplies the reverse inclusion;
hence
\[
 \partial A\subset\partial_{\mathrm{mt}}A\subset\partial A,
 \qquad
 \partial A\cap(A^{(0)}\cup A^{(1)})=\varnothing.
\]
\end{proof}

\begin{proof}[Proof of \cref{thm:main}]
To transfer the preceding conclusions from the normalized problem to
circulation one, we use \cref{prop:turkington} together with the energy
normalization:
\[
 \begin{gathered}
 m=\lambda^{-1},\qquad
 \omega_\lambda=\lambda\one_{A_\lambda},\qquad
 \cE(\one_{A_\lambda})
 =\lambda^{-2}\cE_{\mathrm{kin}}(\omega_\lambda)
 =\max_{\rho\in\mathcal A_m}\cE(\rho),\\
 q_\lambda=\mu_\lambda/\lambda>0,\qquad
 u_\lambda=\mathcal G\one_{A_\lambda}-q_\lambda,\qquad
 A_\lambda=\{u_\lambda>0\}\subset D.
 \end{gathered}
\]
Applying \cref{prop:exhaustion,cor:boundary-density} to this normalized
maximizer gives all the asserted conclusions.
\end{proof}

\begin{remark}\label[remark]{rem:connectedness}
Simple connectedness also follows directly from the equation, without
using boundary regularity or finiteness: by the planar Jordan-curve
criterion, it suffices to show that every smooth Jordan curve in a
component encloses only points of that component. Given a component
$A_j$ and such a curve $\gamma\subset A_j$, simple connectedness of $D$
places its closed interior in $D$, allowing us to set
\[
 \Omega:=\Jint(\gamma),\qquad \overline\Omega\subset D,\qquad
 c_\gamma:=\min_\gamma u>0,\qquad
 \eta=(c_\gamma-u)_+\in H_0^1(\Omega).
\]
The positive minimum on $\gamma$ makes $\eta$ vanish on
$\partial\Omega$, so it is an admissible weak test function. Testing
$-\Delta u=\one_A$ gives
\[
 -\int_\Omega|\nabla\eta|^2\,\mathrm dx
 =\int_\Omega\one_A\eta\,\mathrm dx\geq0
 \quad\Longrightarrow\quad
 u\geq c_\gamma>0\quad\text{on }\overline\Omega.
\]
Continuity extends the inequality throughout $\overline\Omega$;
since this set is connected and meets $A_j$ along $\gamma$, we obtain
\[
 \overline\Omega\subset A
 \quad\Longrightarrow\quad \overline\Omega\subset A_j.
\]
The planar Jordan-curve criterion gives simple connectedness of $A_j$.
\end{remark}


\section{Dirichlet zero modes and convexity}\label{ker:section}
In this section and \cref{an:section}, $D$ is a bounded simply connected
domain with $C^\infty$ boundary; convexity is imposed where stated.
The extension of fixed-mass finiteness and compact-interval bounds to
arbitrary bounded convex domains is proved in \cref{end:section}.
The regularity theorem above applies to every bounded Jordan domain.
We combine the fixed-mass trace inequality with a nodal argument as in
\cite{CabreChanillo}.
Fix a maximizing patch $A$, and write
$\Gamma=\partial A$, $u=\G\one_A-q$, and $g=|\nabla u|$ on $\Gamma$.
\Cref{thm:main} gives $\min_\Gamma g>0$.
Define the single-layer trace and the constrained quadratic form by
\begin{align}
 S_\Gamma\phi&=\left.\G(\phi\,ds|_\Gamma)\right|_\Gamma,
 \qquad X_\Gamma=\left\{\phi\in L^2(\Gamma):\int_\Gamma\phi\ds=0\right\},
 \label{ker:single}\\
 Q_A(\phi)&=\int_\Gamma g\phi^2\ds-
                 \int_\Gamma\phi S_\Gamma\phi\ds\geq0
 \qquad(\phi\in X_\Gamma).\label{ker:Q}
\end{align}
The inequality is \cref{prop:shape} on the full boundary. Set
\[
 K_A=\{\phi\in X_\Gamma:Q_A(\phi)=0\},\qquad
 K_A^0=\{\phi\in X_\Gamma:S_\Gamma\phi=g\phi\}.
\]
Polarization on $X_\Gamma$ gives
\begin{equation}\label{ker:constant}
 \phi\in K_A\quad\Longrightarrow\quad
 S_\Gamma\phi-g\phi=c\quad\hbox{for a constant }c.
\end{equation}
Write $c_A(\phi)=c$; then $K_A^0=\ker c_A$. We call these
the threshold-preserving zero modes, since their first-order
variation leaves the level $q$ fixed; see \eqref{an:shape-derivative}.

\begin{lemma}\label[lemma]{ker:bootstrap}
Every $\phi\in K_A$ is smooth. For $0<\alpha<1$,
\[
 \|\phi\|_{C^{1,\alpha}(\Gamma)}+|c_A(\phi)|
 \leq C_{A,\alpha}\|\phi\|_{L^2(\Gamma)}.
\]
\end{lemma}
\begin{proof}
The logarithmic Green kernel on a fixed smooth curve satisfies
\[
 \begin{aligned}
 \sup_{x\in\Gamma}\|G_D(x,\cdot)\|_{L^2(\Gamma)}&\leq C,\\
 \|G_D(x,\cdot)-G_D(y,\cdot)\|_{L^2(\Gamma)}^2
 &\leq C\delta(1+|\log\delta|)^2,\qquad \delta=|x-y|\ll1.
 \end{aligned}
\]
The second estimate follows by splitting the local arcs at scale $\delta$:
the near part is $O(\delta(1+|\log\delta|)^2)$ and the far part is bounded
by $C\delta^2\int_\delta^1r^{-2}\,dr$. Thus, for $0<\beta<1/2$,
\[
 \|S_\Gamma\phi\|_{C^{0,\beta}(\Gamma)}\leq C_\beta\|\phi\|_2,
\]
\[
 c_A(\phi)=\frac1{\Hh^1(\Gamma)}\int_\Gamma(S_\Gamma\phi-g\phi)\ds,
 \qquad |c_A(\phi)|\leq C\|\phi\|_2.
\]
One-sided Dirichlet regularity gives $g\in C^\infty(\Gamma)$ and
$\min_\Gamma g>0$. Consequently,
\[
 \phi=\frac{S_\Gamma\phi-c_A(\phi)}g,\qquad
 S_\Gamma:C^{k,\beta}(\Gamma)\longrightarrow C^{k+1,\beta}(\Gamma)
 \quad(k\geq0)
\]
imply $\phi\in C^\infty(\Gamma)$. Quantitatively, fix $\beta\in(0,1/2)$
and use
\[
 \|\phi\|_{C^{0,\alpha}}
 \leq C_\alpha\|\phi\|_{C^{1,\beta}}
 \leq C_{A,\alpha}\|\phi\|_2,
 \qquad
 \|\phi\|_{C^{1,\alpha}}\leq
 C_{A,\alpha}(\|\phi\|_{C^{0,\alpha}}+|c_A(\phi)|).
\]
\end{proof}

For $\phi\in K_A^0$, its potential $v=\G(\phi\,ds|_\Gamma)$ satisfies
\begin{equation}\label{ker:weak}
 v\in H^1_0(D),\quad
 \int_D\nabla v\cdot\nabla z\dd=
 \int_\Gamma\frac{vz}{g}\ds\quad(z\in H^1_0(D)),
 \qquad \int_\Gamma\frac{v}{g}\ds=0.
\end{equation}
Conversely, \eqref{ker:weak} gives $\phi=v/g\in K_A^0$.
By \cref{ker:bootstrap}, $v$ is continuous on $\overline D$, smooth up to
each side of $\Gamma$, and harmonic in $D\setminus\Gamma$.

We now use the trace equation to control the nodal domains of its
potential and the zeros of its outer normal derivative.

\begin{lemma}\label[lemma]{ker:nodal}
Every nonzero solution of \eqref{ker:weak} has exactly two nodal domains.
\end{lemma}
\begin{proof}
The trace of $v$ takes both signs by \eqref{ker:weak}. For a nodal domain
$O_i$ of sign $s_i\in\{-1,1\}$ and $\varepsilon>0$, continuity of $v$
on $\overline D$ and its zero boundary values imply that
\[
 K_{i,\varepsilon}:=
 \{x\in\overline O_i:s_iv(x)\geq\varepsilon\}
 \quad\text{is a compact subset of }O_i.
\]
Choose $\chi_{i,\varepsilon}\in C_c^\infty(O_i)$ with
$0\leq\chi_{i,\varepsilon}\leq1$ and
$\chi_{i,\varepsilon}=1$ in a neighborhood of $K_{i,\varepsilon}$.
This cutoff realizes the nodal truncation without assuming any
regularity of $\partial O_i$:
\[
 \begin{gathered}
 v_{i,\varepsilon}:=
 s_i(s_iv-\varepsilon)_+\chi_{i,\varepsilon}
 =s_i(s_iv-\varepsilon)_+\one_{O_i}\in H_0^1(D),\\
 \nabla v_{i,\varepsilon}
 =\one_{O_i\cap\{s_iv>\varepsilon\}}\nabla v.
 \end{gathered}
\]
Indeed, the term containing $\nabla\chi_{i,\varepsilon}$ vanishes
where the positive part is nonzero. For $v_i:=v\one_{O_i}$,
\[
 \begin{gathered}
 \|v_{i,\varepsilon}-v_i\|_{L^2(D)}
 \leq\varepsilon|D|^{1/2},\\
 \|\nabla v_{i,\varepsilon}-\one_{O_i}\nabla v\|_{L^2(D)}^2
 =\int_{O_i\cap\{0<s_iv\leq\varepsilon\}}|\nabla v|^2\dd
 \longrightarrow0.
 \end{gathered}
\]
Thus the truncations converge strongly in $H^1(D)$, and
\[
 v_i\in H_0^1(D),\qquad
 \nabla v_i=\one_{O_i}\nabla v,
\]
so testing \eqref{ker:weak} with $v_i$ gives
\[
 \int_D|\nabla v_i|^2\dd=\int_\Gamma v_i^2/g\ds>0,\qquad
 I_i:=\int_\Gamma v_i/g\ds\ne0.
\]
If three nodal domains existed, put $w=I_2v_1-I_1v_2$. Then $w\ne0$,
$w=0$ on $O_3$, and
\[
 \int_\Gamma w/g\ds=0,\qquad
 \int_D|\nabla w|^2\dd=\int_\Gamma w^2/g\ds.
\]
Polarizing \eqref{eq:weighted-trace} at $w$ gives
\[
 \int_D\nabla w\cdot\nabla z\dd
 =\int_\Gamma w(z-\bar z)/g\ds
 =\int_\Gamma wz/g\ds\qquad(z\in H_0^1(D)).
\]
Here $\bar z=(\int_\Gamma z/g)/(\int_\Gamma1/g)$.

If $w$ vanishes on an open subset of a phase component, harmonic unique
continuation makes it vanish throughout that component. On each adjacent
interface arc its trace, hence the source $w/g$, is zero. The weak equation
then makes $w$ harmonic across the arc, propagating the zero to the next
component. The finitely many Jordan interiors and their connected exterior
have a connected adjacency graph. Since $O_3\setminus\Gamma$ contains an
open set, $w\equiv0$ in $D$, a contradiction.
\end{proof}

\begin{proposition}\label[proposition]{ker:boundaryzeros}
If $0\ne\phi\in K_A^0$ and $v=\G(\phi\,ds|_\Gamma)$, then
$\partial_nv$ has exactly two zeros on $\partial D$. Both are simple
in the boundary arclength parameter. In particular, $\dim K_A^0\leq2$.
\end{proposition}
\begin{proof}
Let $F:\mathbb D\to D$ be conformal. Its smooth boundary extension has
$F'\ne0$ on $\partial\mathbb D$ \cite{Pommerenke}. The function
$w=v\circ F$ is harmonic in an annular collar and vanishes on the circle.
Odd harmonic reflection therefore makes
\[
 b(\theta):=\partial_rw(e^{i\theta})
 =|F'(e^{i\theta})|\,\partial_nv(F(e^{i\theta}))
\]
real analytic. Moreover,
\[
 \int_0^{2\pi}b(\theta)\,d\theta
 =\int_{\partial D}\partial_nv\ds
 =-\int_\Gamma\phi\ds=0.
\]
If $b\equiv0$, zero Cauchy data and the continuation argument in
\cref{ker:nodal} give $v\equiv0$. Thus $b$ takes both signs and has
finitely many zeros. The conformal homeomorphism preserves the two
nodal domains of $v$.

Suppose $b$ has a zero of order $k\geq2$ at $p$ on the circle.
In a local conformal coordinate $\zeta=T(x)$ mapping the disk to the
upper half-plane and $p$ to zero, reflection gives
\[
 w\circ T^{-1}(\zeta)
 =c\,\operatorname{Im}(\zeta^{k+1})+O(|\zeta|^{k+2}),
 \qquad c\in\R\setminus\{0\}.
\]
There are at least three alternating sign sectors. Choose rays
$0<\theta_1<\theta_2<\pi$ in two sectors of the same sign $s$, enclosing
an opposite-sign sector. For small $r>0$, the arcs
\[
 \alpha_i(t)=T^{-1}(te^{i\theta_i}),\qquad0\leq t\leq r,
 \qquad i=1,2,
\]
meet only at $p$ and otherwise lie in the unique sign-$s$ nodal domain $N$.
Choose a simple path $\eta:[0,1]\to N$ joining $\alpha_1(r)$ to
$\alpha_2(r)$, and trim it at
\[
 \begin{split}
 t_2&=\min\{t:\eta(t)\in\alpha_2([0,r])\},\\
 t_1&=\max\{t\leq t_2:\eta(t)\in\alpha_1([0,r])\},\qquad
 \eta(t_i)=\alpha_i(r_i).
 \end{split}
\]
Because $\eta$ avoids $p$, one has $t_1<t_2$ and $r_i>0$. Hence
\[
 J=\alpha_1([0,r_1])\cup\eta([t_1,t_2])\cup\alpha_2([0,r_2])
\]
is a Jordan curve with $J\setminus\{p\}\subset N$ and
$J\cap\partial\mathbb D=\{p\}$. Its interior $\Omega_J$ lies in the disk.
Near $p$, the complementary sector outside $(\theta_1,\theta_2)$ contains
points outside the disk and is therefore on the unbounded side of $J$.
The intervening sign-$-s$ sector lies in $\Omega_J$.

Since $b$ takes both signs, there is a compact boundary subarc $I$ away
from $p$ on which $s b>0$. Uniformly for $\xi\in I$,
\[
 w((1-\delta)\xi)=-\delta b(\xi)+O(\delta^2),
 \qquad\delta\downarrow0.
\]
These collar points have sign $-s$ and lie outside $J$, since
$\dist(I,J)>0$. A path joining them to the sign-$-s$ sector inside $J$
would cross $J\setminus\{p\}\subset N$. This contradicts the uniqueness
of the sign-$-s$ nodal domain. Every zero of $b$ is simple.

If there were at least four zeros, choose cyclically ordered boundary
points $\xi_1,\ldots,\xi_4$ whose inner collars have signs $s,-s,s,-s$.
Join short inward radial segments at $\xi_1,\xi_3$ by a simple path in
$N$, and apply the same first-hit/last-hit trimming. This gives a simple
crosscut $C$ joining $\xi_1$ to $\xi_3$, with its interior in $N$.
Adjoining either circle arc between the endpoints gives a Jordan curve;
therefore $C$ separates the two boundary arcs, and the sign-$-s$ collar
points at $\xi_2,\xi_4$ lie in different components of
$\mathbb D\setminus C$. They cannot lie in one sign-$-s$ nodal domain,
a contradiction. Thus $b$ has exactly two simple zeros. The smooth
boundary parametrization and the positive factor $|F'|$ preserve their
simplicity for $\partial_nv$.

For fixed $\theta_0$, the linear map
\[
 K_A^0\longrightarrow\R^2,\qquad
 \phi\longmapsto\bigl(b(\theta_0),b'(\theta_0)\bigr)
\]
is injective, since a nonzero element of its kernel would give a multiple
boundary zero. Hence $\dim K_A^0\leq2$.
\end{proof}

The nodal information becomes an obstruction to zero modes once it
is combined with the following bilinear Pohozaev identity.

\begin{lemma}\label[lemma]{ker:pohozaev}
Let $U=\G\one_A$, $P=-\partial_nU>0$ on $\partial D$, and let
$\phi\in K_A^0$, $v=\G(\phi\,ds|_\Gamma)$, $b=\partial_nv$.
For every smooth vector field $Z$ on $\overline D$ satisfying
$DZ+(DZ)^T=(\operatorname{div}Z)I$, one has
\begin{equation}\label{ker:bilinear}
 \int_{\partial D}Pb(Z\cdot n)\ds
   =-\int_Av\operatorname{div}Z\dd.
\end{equation}
In particular,
\begin{equation}\label{ker:orthogonality}
 \int_{\partial D}Pb(a+\sigma x)\cdot n\ds=0
 \qquad(a\in\R^2,\ \sigma\in\R).
\end{equation}
\end{lemma}
\begin{proof}
For the bilinear stress vector
\[
 B=(Z\cdot\nabla U)\nabla v+(Z\cdot\nabla v)\nabla U
       -(\nabla U\cdot\nabla v)Z,
\]
the terms involving first derivatives of $Z$ cancel. Since
$-\Delta U=\one_A$, $-\Delta v=\phi\,ds|_\Gamma$,
$\nabla U=-g\nu$ on $\Gamma$, and $v=g\phi$ there,
\[
 \operatorname{div}B
   =-\one_A Z\cdot\nabla v+v(Z\cdot\nu)\,ds|_\Gamma
\]
in distributions. This also follows by integration in the two phases
and the jump condition for the normal derivative of $v$.
Integration over $D$ and then by parts in $A$ gives
\[
 \int_{\partial D}B\cdot n\ds=\int_Av\operatorname{div}Z\dd.
\]
On $\partial D$, both potentials vanish and $B\cdot n=-Pb(Z\cdot n)$,
which proves \eqref{ker:bilinear}. Finally, Green symmetry gives
\[
 \int_Av\dd=\int_\Gamma U\phi\ds=q\int_\Gamma\phi\ds=0.
\]
The field $Z=a+\sigma x$ has constant divergence $2\sigma$, proving
\eqref{ker:orthogonality}.
\end{proof}

\begin{lemma}\label[lemma]{ker:separator}
Let $D$ be smooth and convex, and let $b\in C(\partial D)$ have exactly
two zeros and opposite signs on the complementary boundary arcs.
There is a field $Z(x)=a+\sigma x$ such that $b(Z\cdot n)\geq0$
on $\partial D$, with strict inequality on a nonempty open arc.
\end{lemma}
\begin{proof}
Let $p_1,p_2$ be the zeros. Write $\ell_i$ for the supporting tangent
line at $p_i$ and $F_i=\ell_i\cap\overline D$ for its contact face,
which can be a point or a closed segment.

If $\ell_1=\ell_2$, set $z=(p_1+p_2)/2$. For $Z=x-z$, convexity gives
$Z\cdot n\geq0$, with equality exactly on the common face. Indeed,
any supporting line through $z$ is its unique tangent line by smoothness.
On the complementary boundary arc $b$ has a constant nonzero sign.
Changing the sign of $Z$ proves the assertion.

If the supporting lines are distinct and not parallel, let $z$ be their
intersection. Then $z\notin\overline D$, since two distinct supporting
lines cannot pass through a smooth boundary point. From this exterior
point there are exactly two supporting lines. Thus $(x-z)\cdot n$
vanishes precisely on $F_1\cup F_2$. It is negative at a nearest
boundary point to $z$ and positive at a farthest one, so it has opposite
signs on the two open arcs connecting the faces. The signs of $b$ are
also opposite and constant on these arcs. One overall choice of sign
for $Z=x-z$ gives the required comparison.

If the supporting lines are distinct and parallel, their outer normals
are opposite. Choose a nonzero constant vector $a$ parallel to them.
The zero set of $a\cdot n$ is exactly $F_1\cup F_2$, and its signs on
the two connecting arcs are opposite. The same comparison applies.
\end{proof}

\begin{proposition}\label[proposition]{ker:dirichlet}
In a smooth bounded convex domain, every maximizing patch satisfies
$K_A^0=\{0\}$. Consequently, the map $\phi\mapsto c$ in
\eqref{ker:constant} is injective on $K_A$, and $\dim K_A\leq1$.
\end{proposition}
\begin{proof}
If $0\ne\phi\in K_A^0$, Proposition~\ref{ker:boundaryzeros} supplies
the two sign changes required by Lemma~\ref{ker:separator}.
The resulting field $Z=a+\sigma x$ satisfies
\[
 \int_{\partial D}Pb(Z\cdot n)\ds>0,
\]
because $P>0$. This contradicts \eqref{ker:orthogonality}.
The assertion about the full constrained kernel follows immediately
from \eqref{ker:constant}.
\end{proof}

\section{Analytic reduction and finiteness}\label{an:section}

Throughout this section, $D$ is smooth, bounded, and convex. For
$A\in\M_m$, write
\[
 U_A=\G\one_A,\qquad q_A>0,\qquad
 \{U_A=q_A\}=\partial A,\qquad
 g_A=-\partial_\nu U_A>0\quad\text{on }\partial A.
\]

\begin{lemma}\label[lemma]{an:compactness}
For every compact interval $I\subset(0,|D|)$, the set
\[
 \mathscr M_I=\{(m,\one_A):m\in I,\ A\in\M_m\}
\]
is compact in $I\times L^1(D)$, and $(m,A)\mapsto q_A$ is continuous.
If $(m_j,\one_{A_j})\to(m,\one_A)$, then $\partial A_j$ is eventually
a normal graph over $\partial A$, with graph functions converging to
zero in $C^{1,\alpha}(\partial A)$ for every $0<\alpha<1$.
\end{lemma}

\begin{proof}
Let $\tau=\G1$. Green symmetry and $0\leq\G(\rho+\sigma)\leq2\tau$ give
\[
 \begin{aligned}
 |\cE(\rho)-\cE(\sigma)|
 &=\frac12\left|\int_D(\rho-\sigma)\G(\rho+\sigma)\dd\right|\\
 &\leq\|\tau\|_\infty\|\rho-\sigma\|_1
 \qquad(0\leq\rho,\sigma\leq1).
 \end{aligned}
\]
For $\rho\in\mathcal A_m$, choose
\[
 \sigma=
 \begin{cases}
 (m'/m)\rho,&m'<m,\\[1mm]
 \rho+\dfrac{m'-m}{|D|-m}(1-\rho),&m'>m.
 \end{cases}
 \qquad
 \sigma\in\mathcal A_{m'},\quad
 \|\sigma-\rho\|_1=|m'-m|.
\]
Reversing $m,m'$ yields
\begin{equation}\label{an:value-lipschitz}
 |V(m')-V(m)|\leq\|\tau\|_\infty|m'-m|.
\end{equation}
If $m_j\to m$ and $\one_{A_j}\overset{*}{\rightharpoonup}\rho$, weak-* energy
continuity gives $\rho\in\M_m$. Thus $\rho=\one_A$ and
\[
 \|\one_{A_j}-\one_A\|_1
 =m_j+m-2\int_A\one_{A_j}\dd\longrightarrow0.
\]
This proves compactness. For every finite $p$,
\[
 \|\one_{A_j}-\one_A\|_p
 =\|\one_{A_j}-\one_A\|_1^{1/p}\longrightarrow0,
 \qquad
 \|U_{A_j}-U_A\|_{W^{2,p}}
 \leq C_p\|\one_{A_j}-\one_A\|_p.
\]
Consequently
\begin{equation}\label{an:potential-convergence}
 U_{A_j}\longrightarrow U_A
 \quad\text{in }C^{1,\beta}(\overline D)
 \qquad(0<\beta<1).
\end{equation}

Since $0<q_{A_j}<\|U_{A_j}\|_\infty\leq\|\tau\|_\infty$, every
subsequence has a threshold cluster value $q_*$. Uniform convergence gives
\[
 \one_A=1\ \text{a.e. on }\{U_A>q_*\},\qquad
 \one_A=0\ \text{a.e. on }\{U_A<q_*\}.
\]
If $q_*\neq q_A$, a sufficiently thin strip along the noncritical
level $\partial A$ contradicts one of these identities. Hence
$q_{A_j}\to q_A$.

In disjoint normal tubes around $\partial A$, choose the normal coordinate
$r$ so that
\[
 \partial_r U_A\leq-\tfrac12\min_{\partial A}g_A<0.
\]
Equation \eqref{eq:full-zero-set} bounds $|U_A-q_A|$ away from zero
on the complement of the tubes. Convergence therefore gives exactly one
root $r=h_j(x)$ on every normal fiber and no boundary points elsewhere.
The implicit-function formula and \eqref{an:potential-convergence} yield
\[
 \|h_j\|_{C^1}\to0,\qquad
 \sup_j\|h_j\|_{C^{1,\beta}}<\infty
 \quad\Longrightarrow\quad
 \|h_j\|_{C^{1,\alpha}}\to0\qquad(\alpha<\beta<1).
\]
\end{proof}

Fix $0<\alpha<1$ and $A_0\in\M_{m_0}$, and set
\[
 \Gamma_0=\partial A_0,\qquad q_0=q_{A_0},\qquad
 g_0=-\partial_\nu U_{A_0},\qquad B=C^{1,\alpha}(\Gamma_0).
\]
In disjoint normal tubes, use $x=\gamma(\xi)+r\nu(\xi)$ with
Jacobian $1+\kappa(\xi)r$. For $\|h\|_B$ small, put
\[
 X_h(x)=x+h(x)\nu(x),\qquad
 \partial A_h=X_h(\Gamma_0).
\]
Here $A_h$ is the image of $A_0$ under the $C^{1,\alpha}$
diffeomorphism
\[
 (\xi,r)\longmapsto(\xi,r+\chi(r)h(\xi)),
\]
where $\chi=1$ near $0$ and is supported in the tubes; the map is the
identity elsewhere. Thus $A_h\subset D$ has the topology of $A_0$.
For two small graphs,
\[
 |A_h\mathbin{\triangle}A_k|
 =\int_{\Gamma_0}\left|\int_{k(x)}^{h(x)}(1+\kappa(x)r)\,dr\right|\ds
 \asymp\|h-k\|_{L^1(\Gamma_0)}.
\]
In particular, the chart is injective modulo null sets. Define
\begin{equation}\label{an:equation}
 \mathscr F(h,t,s)
 =\bigl((\G\one_{A_h})\circ X_h-t,\ |A_h|-s\bigr)
 \in Y:=B\times\R
\end{equation}
near $(0,q_0,m_0)$.

\begin{lemma}\label[lemma]{an:analytic-map}
The map $\mathscr F$ is real analytic. Its shape derivative
\begin{equation}\label{an:shape-derivative}
 L\phi:=D_h\mathscr F(0,q_0,m_0)\phi
 =\left(S_{\Gamma_0}\phi-g_0\phi,\ \int_{\Gamma_0}\phi\ds\right)
\end{equation}
is injective and Fredholm of index $-1$.
\end{lemma}

\begin{proof}
Set $p(y)=-|y|^2/4$. Since $-\Delta p=1$, Green's formula gives
\[
 \G\one_{A_h}(x)
 =p(x)+\int_{\partial A_h}
 \bigl[p(y)\partial_{\nu_y}G_D(x,y)
       -G_D(x,y)\partial_{\nu_h}p(y)\bigr]\ds_y
 \qquad(x\in A_h).
\]
Let $S^D$ have kernel $G_D$ and $K^D$ the principal-value kernel
$\partial_{\nu_y}G_D$. The interior double-layer trace is
$-\tfrac12I+K^D$, so
\begin{equation}\label{an:layer-formula}
 (\G\one_{A_h})|_{\partial A_h}
 =\tfrac12p+K^D_{\partial A_h}p
                 -S^D_{\partial A_h}\partial_{\nu_h}p.
\end{equation}

For the fundamental kernel $\Phi(x)=-(2\pi)^{-1}\log|x|$,
\cite[Theorem~3.2]{DRLM} gives real-analytic dependence on
$C^{1,\alpha}$ embeddings in the spaces
\[
 S:C^{0,\alpha}\longrightarrow C^{1,\alpha},\qquad
 K:C^{1,\alpha}\longrightarrow C^{1,\alpha}.
\]
Apply this theorem to each Jordan interior of $\Gamma_0$, whose
exterior is connected. Write $J_h\ds=X_h^*(ds|_{\partial A_h})$ and
$\widehat\nu_h=\nu_h\circ X_h$. For a unit tangent $T_0$ to $\Gamma_0$,
the normal-coordinate convention above gives
\[
 J_h=\sqrt{(1+\kappa h)^2+|\partial_s h|^2},\qquad
 J_h\widehat\nu_h=(1+\kappa h)\nu-(\partial_s h)T_0.
\]
For $\|h\|_B$ small, $J_h$ is bounded away from zero. Products and the
positive square root in the Banach algebra $C^{0,\alpha}(\Gamma_0)$
therefore show that the pulled-back data satisfy
\[
 h\longmapsto p\circ X_h\in C^{1,\alpha},\qquad
 h\longmapsto\widehat\nu_h,\ J_h,\
             (\partial_{\nu_h}p)\circ X_h\in C^{0,\alpha}
\]
real analytically.

The regular Green part $H_D=\Phi-G_D$ is smooth across the diagonal
and harmonic in each variable. Hence
\[
 (\Delta_x+\Delta_y)H_D=0,\qquad
 H_D\in C^\omega(D\times D).
\]
Choose a compact set $K\subset D$ whose interior contains all small
deformed closures. Interior analytic estimates give constants $M,R>0$
such that, for all multiindices $\beta,\gamma\in\mathbb N_0^2$,
\begin{equation}\label{an:green-factorial}
 \sup_{K\times K}|\partial_x^\beta\partial_y^\gamma H_D|
 \leq M R^{-|\beta|-|\gamma|}(|\beta|+|\gamma|)!.
\end{equation}
For distinct components, the same estimate holds for $\Phi(x-y)$ on
the product of their disjoint fixed normal tubes. Thus these bounds
require no regularity of the outer boundary.

To obtain analyticity in operator norm, write
$\Gamma_0=\bigsqcup_{i=1}^N\Gamma_i$ and let $P(x,y)$ denote either
$H_D(x,y)$, or $\Phi(x-y)$ when $i\ne j$. Define
\[
 \begin{aligned}
 T_P^{ij}(h)\eta(x)
 &=\int_{\Gamma_j}P(X_h(x),X_h(y))\eta(y)J_h(y)\ds_y,\\
 K_P^{ij}(h)\xi(x)
 &=\int_{\Gamma_j}\nabla_yP(X_h(x),X_h(y))
       \cdot\widehat\nu_h(y)\xi(y)J_h(y)\ds_y,
 \qquad x\in\Gamma_i.
 \end{aligned}
\]
Set $V_f(x,y)=(f(x)\nu(x),f(y)\nu(y))$. Since $X_h$ is affine in $h$,
\[
 \begin{aligned}
 &D_h^k[P(X_h(x),X_h(y))][f_1,\ldots,f_k]\\
 &\qquad=D^kP(X_h(x),X_h(y))
             [V_{f_1}(x,y),\ldots,V_{f_k}(x,y)].
 \end{aligned}
\]
The bounds \eqref{an:green-factorial}, the H\"older product estimates,
and the formulas for $J_h$ and $J_h\widehat\nu_h$ consequently imply
\begin{equation}\label{an:layer-factorial}
 \begin{aligned}
 &\|D_h^kT_P^{ij}(h)[f_1,\ldots,f_k]\|
       _{\mathcal L(C^{0,\alpha}(\Gamma_j),C^{1,\alpha}(\Gamma_i))}\\
 &\quad+
 \|D_h^kK_P^{ij}(h)[f_1,\ldots,f_k]\|
       _{\mathcal L(C^{1,\alpha}(\Gamma_j),C^{1,\alpha}(\Gamma_i))}
 \leq C_0C_1^k k!\prod_{\ell=1}^k\|f_\ell\|_B.
 \end{aligned}
\end{equation}
The constants are uniform on a smaller ball about $0$ in $B$.
Indeed, differentiating once in the target variable and estimating
its H\"older seminorm require at most three additional derivatives
of $P$; the resulting factor $(k+3)!$ is absorbed into $C_0C_1^k k!$.
For either operator $\mathcal T=T_P^{ij}$ or $K_P^{ij}$, Taylor's formula
in its indicated operator norm now gives
\[
 \left\|\mathcal T(h)-\sum_{k=0}^n\frac1{k!}
          D_h^k\mathcal T(0)[h,\ldots,h]\right\|
 \leq C_0(C_1\|h\|_B)^{n+1}\longrightarrow0
 \qquad(C_1\|h\|_B<1).
\]
Together with the singular diagonal blocks supplied by \cite{DRLM},
this proves analyticity of the full pulled-back layer operators.
Equation \eqref{an:layer-formula}, together with
\[
 |A_h|=m_0+\int_{\Gamma_0}
                  \left(h+\frac{\kappa h^2}{2}\right)\ds,
\]
proves analyticity of $\mathscr F$.

For $h=\varepsilon\phi$, write $U_\varepsilon=\G\one_{A_{\varepsilon\phi}}$.
The signed-strip formula yields, at each $x\in\Gamma_0$,
\[
 \frac{U_\varepsilon(x)-U_0(x)}{\varepsilon}
 \longrightarrow\int_{\Gamma_0}G_D(x,y)\phi(y)\ds_y.
\]
Indeed, away from an arc of length $2\delta$ about $x$, ordinary transport
applies; the arc contributes at most
$C\delta(1+|\log\delta|)$ to the difference quotient for small
$|\varepsilon|$. First let $\varepsilon\to0$, then $\delta\to0$.
Also $U_\varepsilon\to U_0$ in $C^1(\overline D)$ by Dirichlet estimates.
Therefore
\[
 \begin{aligned}
 \frac{U_\varepsilon(X_{\varepsilon\phi}(x))-U_0(x)}{\varepsilon}
 &=
 \frac{U_\varepsilon(x)-U_0(x)}{\varepsilon}\\
 &\quad+\int_0^1\nabla U_\varepsilon(x+\theta\varepsilon\phi(x)\nu(x))
                  \cdot\phi(x)\nu(x)\,d\theta\\
 &\longrightarrow S_{\Gamma_0}\phi(x)-g_0(x)\phi(x).
 \end{aligned}
\]
Analyticity into $B$ identifies this pointwise formula with the
Banach-space derivative, proving \eqref{an:shape-derivative}.

One-sided elliptic regularity gives $g_0\in C^\infty(\Gamma_0)$,
and $\min_{\Gamma_0}g_0>0$. The operator
\[
 L_0\phi=\left(-g_0\phi,\int_{\Gamma_0}\phi\ds\right)
\]
is an isomorphism from $B$ onto
\[
 \left\{(f,b)\in Y:\ b+\int_{\Gamma_0}\frac f{g_0}\ds=0\right\}.
\]
Thus $\operatorname{ind}L_0=-1$. Since
\[
 B\hookrightarrow C^{0,\alpha}(\Gamma_0)
 \quad\text{compactly},\qquad
 S_{\Gamma_0}:C^{0,\alpha}\longrightarrow B
 \quad\text{boundedly},
\]
$L-L_0$ is compact and $\operatorname{ind}L=-1$. Finally,
\[
 L\phi=0
 \quad\Longrightarrow\quad
 S_{\Gamma_0}\phi=g_0\phi,\quad\int_{\Gamma_0}\phi\ds=0
 \quad\Longrightarrow\quad
 \phi\in K_{A_0}^0=\{0\}
\]
by \cref{ker:dirichlet}.
\end{proof}

\begin{proposition}[Scalar analytic reduction]\label[proposition]{an:scalar}
There are open intervals $J\ni q_0$, $I_0\ni m_0$, a neighborhood
$\mathcal U\ni0$ in $B$, and real-analytic maps
$H:J\times I_0\to B$, $r:J\times I_0\to\R$, with
$H(q_0,m_0)=0$, such that
\begin{equation}\label{an:scalar-equation}
 \mathscr F(h,t,s)=0
 \quad\Longleftrightarrow\quad
 h=H(t,s),\quad r(t,s)=0
 \qquad((h,t,s)\in\mathcal U\times J\times I_0).
\end{equation}
At fixed mass and threshold, there is at most one maximizing patch
in this chart.
\end{proposition}

\begin{proof}
The range $Y_0=L(B)$ is closed of codimension one. Choose
\[
 \ell\in Y^*,\quad\ker\ell=Y_0,\qquad
 e\in Y,\quad\ell(e)=1,\qquad Py=y-\ell(y)e.
\]
Then $D_h(P\mathscr F)=L:B\to Y_0$ is an isomorphism.
The analytic implicit-function theorem gives $h=H(t,s)$ as the unique
solution of $P\mathscr F=0$ in product neighborhoods. Shrink them so
that $H(J\times I_0)\subset\mathcal U$, $J\subset(0,\infty)$, and
$I_0\subset(0,|D|)$. Setting
\[
 r(t,s)=\ell\mathscr F(H(t,s),t,s)
 \quad\text{gives}\quad
 \mathscr F(H(t,s),t,s)=r(t,s)e,
\]
which proves \eqref{an:scalar-equation}. By \cref{an:compactness},
all maximizing pairs sufficiently close to $(m_0,\one_{A_0})$ in
mass and $L^1$ lie in this chart, with thresholds in $J$.
Injectivity of the graph chart gives the last assertion.
\end{proof}

\begin{lemma}\label[lemma]{an:constant-energy}
If $r(t,m_0)$ vanishes identically near $q_0$, then, on a smaller
open interval $J_1\ni q_0$, the patches $A(t):=A_{H(t,m_0)}$ are
maximizers with thresholds $t$. They are distinct for distinct $t$
and continuous in $L^1(D)$.
\end{lemma}

\begin{proof}
The full equation gives
\[
 |A(t)|=m_0,\qquad U_{A(t)}=t\quad\text{on }\partial A(t).
\]
For their normal velocity $V_t$,
\[
 \frac d{dt}\cE(\one_{A(t)})
 =\int_{\partial A(t)}U_{A(t)}V_t\ds
 =t\int_{\partial A(t)}V_t\ds
 =t\frac d{dt}|A(t)|=0.
\]
To justify the first equality, put
$f=\one_{A(t+\delta)}-\one_{A(t)}$. Since $H$ is $C^1$,
\[
 |A(t+\delta)\mathbin{\triangle}A(t)|=O(|\delta|),\qquad
 \|f\|_{\dot H^{-1}_D}^2
 \leq C_p\|f\|_p^2=O(|\delta|^{2/p})=o(|\delta|)
 \quad(1<p<2).
\]
The exact energy expansion reduces the derivative to transport of
the continuous potential $U_{A(t)}$. As $A(q_0)=A_0$,
\[
 \cE(\one_{A(t)})=V(m_0),\qquad A(t)\in\M_{m_0}.
\]
Its unique threshold is $t$: if $q_{A(t)}\neq t$, continuity
of $U_{A(t)}$ and positive-area neighborhoods on both sides of its
$C^{1,\alpha}$ boundary contradict
$\one_{A(t)}=\one_{\{U_{A(t)}>q_{A(t)}\}}$ a.e.
Threshold uniqueness and the signed-strip estimate now give
distinctness and $L^1$ continuity.
\end{proof}

\begin{proposition}\label[proposition]{an:finiteness}
For every $m\in(0,|D|)$, the set $\M_m$ is finite.
\end{proposition}

\begin{proof}
Suppose the compact metric space $(\M_m,d)$,
$d(A,B)=\|\one_A-\one_B\|_1$, is infinite. Its derived set
$\M_m'$ is nonempty and compact. Choose
\[
 A_0\in\M_m',\qquad q_0=q_{A_0}=\max_{A\in\M_m'}q_A,
 \qquad A_j\to A_0,
\]
with $A_j\in\M_m\setminus\{A_0\}$ pairwise distinct. By the local
reduction, for large $j$,
\[
 A_j=A_{H(q_{A_j},m)},\qquad r(q_{A_j},m)=0,\qquad q_{A_j}\to q_0.
\]
The thresholds are pairwise distinct, since fixed mass and threshold
determine at most one shape in this chart. The one-variable identity
theorem gives $r(t,m)\equiv0$ near $q_0$.
By \cref{an:constant-energy}, there is an open interval $J_1\ni q_0$
of distinct maximizers $A_{H(t,m)}$ with thresholds $t$. Their
$L^1$ continuity implies
\[
 A_{H(t,m)}\in\M_m'\qquad(t\in J_1).
\]
Taking $t>q_0$ contradicts the choice of $A_0$.
\end{proof}

\begin{corollary}\label[corollary]{an:uniform-interior}
For every compact interval $I\subset(0,|D|)$ there is $N_I<\infty$ such that
\[
 \#\M_m\leq N_I\qquad(m\in I).
\]
\end{corollary}

\begin{proof}
At $(m_0,A_0)\in\mathscr M_I$, fixed-mass finiteness and
\cref{an:constant-energy} imply $r(\cdot,m_0)\not\equiv0$.
Its zero at $q_0$ has finite order $k\geq1$. Weierstrass preparation gives
\[
 r(t,s)=a(t,s)\left[(t-q_0)^k+
             \sum_{j=0}^{k-1}a_j(s)(t-q_0)^j\right],
 \qquad a\neq0,\quad a_j(m_0)=0.
\]
Each mass fiber has at most $k$ roots, hence at most $k$ maximizing
patches in this chart. By \cref{an:compactness}, the smaller chart on which
this preparation holds contains
a relatively open neighborhood of $(m_0,A_0)$ in $\mathscr M_I$.
A finite subcover, with orders $k_1,\ldots,k_L$, yields
\[
 \#\M_m\leq\sum_{\ell=1}^Lk_\ell=:N_I\qquad(m\in I).
\]
\end{proof}

\section{Uniform counting}\label{end:section}

We first prove \cref{rem:convex-jordan-finiteness}, allowing any bounded
convex fluid domain. The boundary argument is transferred to smooth
convex subdomains, while the shape analysis remains confined to the
interior. We then impose a smooth boundary to control the two mass
endpoints and prove \cref{thm:convex-finiteness}.

\begin{lemma}\label[lemma]{end:rough-kernel}
Let $D\subset\R^2$ be a bounded convex domain and $A\in\M_m$, where
$0<m<|D|$. With the notation of \eqref{ker:single}--\eqref{ker:constant},
\[
 K_A^0=\{0\},\qquad \dim K_A\leq1.
\]
\end{lemma}

\begin{proof}
Suppose $0\ne\phi\in K_A^0$, and set
\[
 \Gamma=\partial A,\qquad U=\G\one_A,\qquad
 v=\G(\phi\,ds|_\Gamma),\qquad g=|\nabla U|\big|_\Gamma.
\]
By \cref{thm:main}, $\Gamma$ is a finite union of analytic curves
contained in $D$, and $g>0$. The proof of \cref{ker:bootstrap} uses only
layer estimates on these interior curves, so $\phi\in C^\infty(\Gamma)$.
The interior trace theorem and Poincar\'e's inequality give
\[
 \left|\int_\Gamma\phi z\ds\right|
 \leq C_\Gamma\|\phi\|_{L^2(\Gamma)}\|\nabla z\|_{L^2(D)}
 \qquad(z\in H_0^1(D)),
\]
so the curve source belongs to $H^{-1}(D)$ and defines a unique weak
Dirichlet potential. To verify its continuity at the outer boundary,
consider the Newtonian single-layer potential
\[
 N_\phi(x)=\int_\Gamma\Phi(x-y)\phi(y)\ds_y.
\]
Single-layer regularity on the smooth interior curves gives a
continuous potential with $C^{1,\alpha}$ extensions from either side
of each curve \cite[Section~2]{DRLM}. Since $N_\phi$ is smooth away
from $\Gamma$, its matching traces and bounded one-sided gradients imply
\[
 N_\phi\in H^1(D)\cap C(\overline D),\qquad
 -\Delta N_\phi=\phi\,ds|_\Gamma\quad\text{in }\mathcal D'(D).
\]

Identify $\R^2$ with $\mathbb C$, write $\mathbb T=\partial\mathbb D$,
and let $F:\mathbb D\to D$ be conformal. Its continuous extension is a
homeomorphism of the closures. Let $\mathcal P$ denote the Poisson
extension on $\mathbb D$, and set
\[
 W_\phi=N_\phi\circ F,\qquad
 h_\phi=\mathcal P[W_\phi|_{\mathbb T}].
\]
Conformal invariance of Dirichlet energy and the Dirichlet principle
on the disk yield
\[
 \begin{gathered}
 W_\phi,h_\phi\in H^1(\mathbb D)\cap C(\overline{\mathbb D}),\qquad
 W_\phi-h_\phi\in H_0^1(\mathbb D),\\
 \int_{\mathbb D}|\nabla h_\phi|^2\dd
 \leq\int_{\mathbb D}|\nabla W_\phi|^2\dd
 =\int_D|\nabla N_\phi|^2\dd.
 \end{gathered}
\]
Choose $f_k\in C_c^\infty(\mathbb D)$ converging to
$W_\phi-h_\phi$ in $H^1(\mathbb D)$. Then
\[
 \begin{gathered}
 f_k\circ F^{-1}\in C_c^\infty(D),\\
 \|\nabla(f_k\circ F^{-1}-f_\ell\circ F^{-1})\|_{L^2(D)}
 =\|\nabla(f_k-f_\ell)\|_{L^2(\mathbb D)}\longrightarrow0.
 \end{gathered}
\]
Poincar\'e's inequality makes these pullbacks Cauchy in $H_0^1(D)$;
local convergence identifies their limit as
\[
 \widetilde v=N_\phi-h_\phi\circ F^{-1}
 \in H_0^1(D)\cap C(\overline D),\qquad
 \widetilde v|_{\partial D}=0.
\]
The harmonicity of $h_\phi\circ F^{-1}$ gives
$-\Delta\widetilde v=\phi\,ds|_\Gamma$. Weak uniqueness therefore
identifies $\widetilde v=v$, and the defining equation for $K_A^0$ gives
\[
 v\in H_0^1(D)\cap C(\overline D),\qquad
 v|_{\partial D}=0,\qquad v=g\phi\quad\text{on }\Gamma.
\]
The truncation and continuation proof of \cref{ker:nodal} therefore
applies without a smooth outer boundary: $v$ has exactly two nodal
domains. Now set
\[
 W=U\circ F,\qquad w=v\circ F.
\]
Both functions are harmonic in an annulus $r_0<|z|<1$ and vanish on
$\mathbb T$. Reflection across the circle makes their radial
derivatives real analytic there. The Hopf lemma on the disk collar and
the zero total source give
\begin{equation}\label{end:circle-data}
 P_*(\theta):=-\partial_rW(e^{i\theta})>0,\qquad
 \beta(\theta):=\partial_rw(e^{i\theta}),\qquad
 \int_0^{2\pi}\beta(\theta)\,d\theta
 =-\int_\Gamma\phi\ds=0.
\end{equation}
The disk argument in the proof of \cref{ker:boundaryzeros} now applies
to $w$: $\beta$ has exactly two simple zeros. This conclusion concerns
the radial derivative on the circle and does not require a pointwise
normal derivative on $\partial D$.

Convexity supplies the uniform estimate needed to pass the stress
identity to the limiting boundary. Indeed, the conformal convexity
criterion and its Herglotz representation give a probability measure
$\mu$ on $\mathbb T$ such that
\[
 \operatorname{Re}\left(1+\frac{zF''(z)}{F'(z)}\right)>0,\qquad
 \log\frac{F'(z)}{F'(0)}
 =-2\int_{\mathbb T}\log(1-\overline\xi z)\,d\mu(\xi).
\]
Taking real parts yields
\begin{equation}\label{end:convex-distortion}
 |F'(z)|\geq\frac{|F'(0)|}{(1+|z|)^2},\qquad
 \left\|\frac1{F'}\right\|_{L^\infty(\mathbb D)}
 \leq\frac4{|F'(0)|}.
\end{equation}
The same convexity criterion shows that
\[
 D_r:=F(r\mathbb D)\quad(0<r<1)
\]
is a smooth convex domain, and $D_r\uparrow D$. Since $\partial D$ is
rectifiable, $F'$ belongs to the Hardy class with exponent one
\cite{Pommerenke}. In particular, its radial boundary values satisfy
\[
 \begin{gathered}
 \sup_{0<r<1}\int_0^{2\pi}|F'(re^{i\theta})|\,d\theta<\infty,\\
 F'(re^{i\cdot})\longrightarrow F'(e^{i\cdot})
 \quad\text{in }L^1(0,2\pi)\text{ and a.e. as }r\uparrow1.
 \end{gathered}
\]
The boundary parametrization $\gamma(\theta)=F(e^{i\theta})$ is
absolutely continuous, with
$\gamma'(\theta)=ie^{i\theta}F'(e^{i\theta})$ a.e. Together with
\eqref{end:convex-distortion}, this gives, for almost every $\theta$
and $\zeta=e^{i\theta}$,
\begin{equation}\label{end:boundary-parametrization}
 0<|F'(\zeta)|<\infty,\qquad
 ds=|F'(\zeta)|\,d\theta,\qquad
 n(F(\zeta))=\frac{\zeta F'(\zeta)}{|F'(\zeta)|}.
\end{equation}
Here $F'(\zeta)$ denotes the radial limit and $n$ the almost-everywhere
outward normal of the convex domain. More explicitly, the arclength
change of variables reads
\begin{equation}\label{end:boundary-change-of-variables}
 \int_{\partial D}f(x)\ds_x
 =\int_0^{2\pi}f(F(e^{i\theta}))|F'(e^{i\theta})|\,d\theta
 \qquad(f\geq0\text{ measurable}).
\end{equation}
Since the weight is positive and finite a.e., a property holds for
arclength-a.e. point of $\partial D$ if and only if its pullback holds
for a.e. $\theta$.

For $Z(x)=a+\sigma x$, $a\in\R^2$, $\sigma\in\R$, introduce the
bilinear stress vector
\[
 \mathcal B_Z=(Z\cdot\nabla U)\nabla v+(Z\cdot\nabla v)\nabla U
              -(\nabla U\cdot\nabla v)Z.
\]
The distributional calculation in \cref{ker:pohozaev} is confined to
the interior and gives
\[
 \operatorname{div}\mathcal B_Z
 =-\one_A Z\cdot\nabla v+v(Z\cdot\nu)\,ds|_\Gamma.
\]
For $r$ sufficiently close to one, $\overline A\subset D_r$; hence
\begin{equation}\label{end:exhaustion-stress}
 \int_{\partial D_r}\mathcal B_Z\cdot n_r\ds
 =2\sigma\int_Av\dd
 =2\sigma q_A\int_\Gamma\phi\ds=0.
\end{equation}
Since $U$ and $v$ need not vanish on $\partial D_r$, we retain all
terms in the stress vector. With
\[
 C_Z(r,\theta):=
 \frac{a+\sigma F(re^{i\theta})}{e^{i\theta}F'(re^{i\theta})},
\]
the normal and tangential components in
\eqref{end:exhaustion-stress} become
\[
 \begin{aligned}
 0={}&\int_0^{2\pi}\operatorname{Re}C_Z
       \bigl(rW_rw_r-r^{-1}W_\theta w_\theta\bigr)\,d\theta\\
 &+\int_0^{2\pi}\operatorname{Im}C_Z
       \bigl(W_\theta w_r+W_rw_\theta\bigr)\,d\theta.
 \end{aligned}
\]
All derivatives here are evaluated at $re^{i\theta}$. Reflection gives,
uniformly in $\theta$,
\[
 W_\theta,w_\theta=O(1-r),\qquad
 W_r=-P_*+O(1-r),\qquad w_r=\beta+O(1-r).
\]
The distortion estimate \eqref{end:convex-distortion} gives
\[
 \sup_{0<r<1}\|C_Z(r,\cdot)\|_{L^\infty(0,2\pi)}
 \leq\frac{4}{|F'(0)|}
       \left(|a|+|\sigma|\sup_{x\in D}|x|\right).
\]
Dominated convergence therefore yields
the boundary identity
\begin{equation}\label{end:rough-orthogonality}
 \int_0^{2\pi}P_*(\theta)\beta(\theta)c_Z(\theta)\,d\theta=0,
 \qquad
 c_Z(\theta):=\operatorname{Re}
 \frac{a+\sigma F(e^{i\theta})}{e^{i\theta}F'(e^{i\theta})}.
\end{equation}

To obtain a field with the required sign, apply \cref{ker:separator}
on $D_r$ to the boundary function
$b_r(F(re^{i\theta}))=\beta(\theta)$. After a positive rescaling, the
resulting field $Z_r(x)=a_r+\sigma_rx$ satisfies
\[
 |a_r|^2+\sigma_r^2=1,\qquad
 \beta(\theta)\operatorname{Re}C_{Z_r}(r,\theta)\geq0.
\]
Choose $r_j\uparrow1$ so that $(a_{r_j},\sigma_{r_j})\to(a,\sigma)$.
Then
\[
 |a|^2+\sigma^2=1,\qquad
 \beta c_Z\geq0\quad\text{a.e. on }\mathbb T.
\]
Since $P_*>0$, \eqref{end:rough-orthogonality} implies
$\beta c_Z=0$ a.e.; the two-zero property of $\beta$ gives $c_Z=0$
a.e. By \eqref{end:boundary-parametrization} and
\eqref{end:boundary-change-of-variables},
\[
 (Z\cdot n)\circ F=|F'|c_Z=0
 \quad\Longrightarrow\quad Z\cdot n=0
 \quad\text{for arclength-a.e. }x\in\partial D.
\]
A bounded convex domain is Lipschitz. The divergence theorem thus gives
\[
 0=\int_{\partial D}Z\cdot n\ds=2\sigma|D|,
 \qquad
 0=\int_{\partial D}(a\cdot x)(a\cdot n)\ds=|a|^2|D|,
\]
where the second equality uses $\sigma=0$. This contradicts
$|a|^2+\sigma^2=1$, proving $K_A^0=\{0\}$. The bound on $\dim K_A$
follows from \eqref{ker:constant}.
\end{proof}

The remaining extension concerns compactness and the analytic shape
map. Uniform convergence of the potentials up to the wall, together
with interior convergence of their gradients, is sufficient for both.

\begin{proof}[Proof of Remark~\ref{rem:convex-jordan-finiteness}]
Let $D$ be bounded and convex. Comparison with the Green function of a
containing disk gives
\[
 0\leq G_D(x,y)\leq C_D\bigl(1+|\log|x-y||\bigr).
\]
H\"older's inequality and interior elliptic estimates imply, for
$1<p<\infty$ and $\overline{D'}\subset D$,
\begin{equation}\label{end:rough-elliptic}
 \|\G f\|_{C(\overline D)}\leq C_{D,p}\|f\|_{L^p(D)},\qquad
 \|\G f\|_{W^{2,p}(D')}\leq C_{D',D,p}\|f\|_{L^p(D)}.
\end{equation}
For general $f\in L^p(D)$, the continuous representative follows by
approximating $f$ by bounded functions and using the first estimate.

Fix a compact interval $I\subset(0,|D|)$. The value estimate
\eqref{an:value-lipschitz} and weak-star energy continuity are valid
on $D$. As in \cref{an:compactness}, they give compactness of
$\mathscr M_I$ in $I\times L^1(D)$. For a convergent sequence
$(m_j,\one_{A_j})\to(m,\one_A)$, \eqref{end:rough-elliptic} gives
\begin{equation}\label{end:rough-potential-convergence}
 U_{A_j}\longrightarrow U_A
 \quad\text{in }C(\overline D)\text{ and in }C^{1,\alpha}_{\mathrm{loc}}(D)
 \qquad(0<\alpha<1).
\end{equation}
Every cluster value $q_*$ of $q_{A_j}$ consequently satisfies
\[
 \one_A=1\ \text{a.e. on }\{U_A>q_*\},\qquad
 \one_A=0\ \text{a.e. on }\{U_A<q_*\}.
\]
If $q_*\ne q_A$, a thin strip on one side of the noncritical level
$\partial A$ contradicts these identities. Hence $q_{A_j}\to q_A$,
and compactness yields
\begin{equation}\label{end:interior-patches}
 q_I:=\min_{(m,A)\in\mathscr M_I}q_A>0,\qquad
 \overline A\subset K_I:=\{x\in\overline D:\tau(x)\geq q_I\}\subset D,
 \qquad \tau=\G1.
\end{equation}
The set $K_I$ is compact. Since $\{U_A=q_A\}=\partial A$, the
potential has a nonzero gap from $q_A$ outside fixed interior normal
tubes. Equations \eqref{end:rough-potential-convergence} and
\eqref{end:interior-patches} therefore give the same normal-graph
convergence as \cref{an:compactness}.

All graphs in a shape chart have closures in a compact subset of $D$.
The regular Green part satisfies $H_D\in C^\omega(D\times D)$, so the
layer-potential proof of analyticity in \cref{an:analytic-map} applies
unchanged. In the derivative calculation, local $C^1$ convergence
from \eqref{end:rough-potential-convergence} replaces convergence on
$\overline D$. Thus
\[
 D_h\mathscr F(0,q_A,m)\phi
 =\left(S_\Gamma\phi-g\phi,\int_\Gamma\phi\ds\right),\qquad
 \operatorname{ind}D_h\mathscr F=-1,\qquad
 \ker D_h\mathscr F=K_A^0=\{0\},
\]
where the last equality is \cref{end:rough-kernel}.
The scalar reduction \eqref{an:scalar-equation} and the constant-energy
argument of \cref{an:constant-energy} are consequently valid on $D$.

These are precisely the compactness and analytic-reduction properties
used in the proof of \cref{an:finiteness}. Its maximum-threshold
argument on the derived set therefore gives $\#\M_m<\infty$ for each
$m\in(0,|D|)$.

For the compact-interval bound, the scalar residual at
$(m_0,A_0)\in\mathscr M_I$ has finite vanishing order
\[
 k_{m_0,A_0}:=\operatorname{ord}_{t=q_{A_0}}r(t,m_0)<\infty;
\]
otherwise \cref{an:constant-energy} would give infinitely many
maximizers at $m_0$. Weierstrass preparation, as in
\cref{an:uniform-interior}, bounds the number of maximizing patches
at each nearby mass in this chart by $k_{m_0,A_0}$. Normal-graph
compactness supplies finitely many such charts covering $\mathscr M_I$,
with orders $k_1,\ldots,k_L$; hence
\[
 1\leq\#\M_m\leq\sum_{\ell=1}^Lk_\ell=:N_{D,I}
 \qquad(m\in I).
\]
\end{proof}

For the rest of this section, $D$ has $C^\infty$ boundary, and convexity
is imposed where stated. We combine the compact-interval count with
uniqueness near $m=0$ and $m=|D|$. The full-mass argument uses a fixed
normal collar and a positive lower bound for the normal derivative of
the torsion function on the entire fluid boundary.

\begin{proposition}\label[proposition]{end:full}
Let $D$ be smooth, bounded, and simply connected. There are
$\delta_D,c_D>0$ such that, for $|D|-\delta_D<m<|D|$,
$\M_m=\{\one_A\}$ and
\[
 Q_A(\phi)\geq c_D\|\phi\|_{L^2(\partial A)}^2
 \qquad\bigl(\phi\in L^2(\partial A)\bigr).
\]
\end{proposition}

\begin{proof}
Set $\tau=\G1$, $\eps=|D|-m$, and fix $p>2$. Uniformly over
$A\in\M_m$, elliptic estimates give
\[
 0\leq U_A\leq\tau,\qquad
 \|U_A-\tau\|_{C^1(\overline D)}
 \leq C_{D,p}\|\one_A-1\|_p
 =C_{D,p}\eps^{1/p},\qquad \|\nabla U_A\|_\infty\leq M.
\]
Moreover, for each fixed $\theta>0$,
\[
 q_A\geq\theta
 \ \Longrightarrow\ A\cap\{\tau\leq\theta\}=\varnothing
 \ \Longrightarrow\ \eps\geq|\{\tau\leq\theta\}|>0.
\]
Hence $q_A\to0$ uniformly as $\eps\downarrow0$.

Use inward normal coordinates $\Psi(s,r)$ in a fixed collar
$0\leq r\leq r_0$, with $s$ arclength on $\partial D$ and
$J(s,r)=|\partial_s\Psi(s,r)|$. The Hopf lemma allows fixed constants
$a_0,J_-,J_+>0$ such that
\[
 \partial_r(\tau\circ\Psi)\geq2a_0,\qquad
 J_-\leq J\leq J_+.
\]
Choose $r_0$ so that the complementary compact core is nonempty.
Since $\min_{\mathrm{core}}\tau>0$, for sufficiently small $\eps$,
\[
 \partial_r(U_A\circ\Psi)\geq a_0\quad\text{in the collar},
 \qquad U_A>q_A\quad\text{on the core}.
\]
Each normal fiber therefore has a unique crossing $r=a(s)\in(0,r_0)$,
and
\[
 D\setminus A=\{\Psi(s,r):0<r\leq a(s)\}.
\]
Write $L=|\partial D|$. The level equation and the mass constraint give
\[
 a_0a(s)\leq q_A\leq Ma(s),\qquad
 \eps=\int_{\partial D}\!\int_0^{a(s)}J(s,r)\,dr\ds
 \geq\frac{J_-L}{M}q_A,
\]
and consequently
\[
 \|a\|_\infty\leq C_w\eps,
 \qquad C_w=\frac{M}{a_0J_-L}.
\]
Differentiating $(U_A\circ\Psi)(s,a(s))=q_A$ yields
\[
 \begin{gathered}
 |a'|\leq\frac{J_+M}{a_0},\qquad
 ds_{\partial A}=j_a(s)\ds,\\
 j_a:=\sqrt{J(s,a)^2+|a'|^2}
 \leq L_0:=J_+\sqrt{1+(M/a_0)^2}.
 \end{gathered}
\]

For $A,B\in\M_m$ with graphs $a,b$, set
$h=b-a$, $f=\one_B-\one_A$, and $\delta=C_w\eps<r_0$.
On the strip between the graphs,
\[
 u_Af\leq0,\qquad
 |u_A(\Psi(s,r))|\geq a_0|r-a(s)|.
\]
Fiberwise integration therefore gives
\begin{equation}\label{end:linear}
 -2\int_Du_Af\dd\geq a_0J_-\|h\|_{L^2(\partial D)}^2.
\end{equation}
For $z\in H_0^1(D)$, its pullback $\widetilde z=z\circ\Psi$ satisfies
\[
 |\widetilde z(s,r)|^2
 \leq\delta\int_0^\delta|\partial_r\widetilde z(s,t)|^2\,dt
 \qquad(0\leq r\leq\delta).
\]
This holds first for smooth zero-trace functions and then by Sobolev
slicing. Cauchy--Schwarz gives
\[
 \begin{aligned}
 \left|\int_Dfz\dd\right|
 &\leq J_+\sqrt\delta\,\|h\|_{L^2(\partial D)}
 \left(\int_{\partial D}\!\int_0^\delta
 |\partial_r\widetilde z|^2\,dr\ds\right)^{1/2}\\
 &\leq\frac{J_+\sqrt\delta}{\sqrt{J_-}}
 \|h\|_{L^2(\partial D)}\|\nabla z\|_2.
 \end{aligned}
\]
Thus
\begin{equation}\label{end:quadratic}
 \|f\|_{\dot H^{-1}_D}^2
 \leq\frac{J_+^2}{J_-}\delta\|h\|_{L^2(\partial D)}^2.
\end{equation}
Since $|A|=|B|$ and $\cE(\one_A)=\cE(\one_B)$,
\eqref{eq:energy-expansion} yields
\[
 -2\int_Du_Af\dd=\|f\|_{\dot H^{-1}_D}^2,\qquad
 \left(a_0J_--\frac{J_+^2}{J_-}\delta\right)
 \|h\|_{L^2(\partial D)}^2\leq0.
\]
Hence $A=B$ for sufficiently small $\delta$.

The same normal-ray estimate gives, with $\Gamma=\partial A$,
\[
 \int_\Gamma z^2\ds
 \leq\frac{L_0}{J_-}\delta\int_D|\nabla z|^2\dd,
 \qquad
 \|\phi\,ds|_\Gamma\|_{\dot H^{-1}_D}^2
 \leq\frac{L_0}{J_-}\delta\|\phi\|_{L^2(\Gamma)}^2.
\]
Since $g\geq a_0$ on $\Gamma$,
\[
 Q_A(\phi)\geq
 \left(a_0-\frac{L_0}{J_-}\delta\right)\|\phi\|_{L^2(\Gamma)}^2.
\]
Choose $\delta_D>0$ below the preliminary mass-deficit bound and so that
\[
 C_w\delta_D\leq
 \min\left\{\frac{r_0}{2},
 \frac{a_0J_-^2}{2J_+^2},\frac{a_0J_-}{2L_0}\right\}.
\]
Both conclusions follow with $c_D=a_0/2$.
\end{proof}

\begin{proposition}\label[proposition]{end:small}
If $D$ is smooth, bounded, and convex, there is $m_D>0$ such that
$\#\M_m=1$ for $0<m<m_D$.
\end{proposition}

\begin{proof}
We verify the concentration, single-core, and strength hypotheses of
\cite[(1.12)--(1.13), Theorems~1.2--1.3]{CaoGuoPengYan}.
For $m_j\downarrow0$ and $A_j\in\M_{m_j}$, put
\[
 \lambda_j=m_j^{-1},\qquad
 \omega_j=\lambda_j\one_{A_j},\qquad
 \psi_j=\lambda_jU_{A_j},\qquad \mu_j=\lambda_jq_{A_j}.
\]
The scaling
\[
 \rho\longmapsto\lambda_j\rho:\mathcal A_{m_j}
 \longrightarrow\{0\leq\omega\leq\lambda_j:\int_D\omega=1\},
 \qquad \cE(\lambda_j\rho)=\lambda_j^2\cE(\rho)
\]
is bijective, so $\omega_j$ is an absolute maximizer in the
circulation-one class and
\[
 -\Delta\psi_j=\lambda_j\one_{A_j},\quad
 \psi_j|_{\partial D}=0,\quad
 A_j=\{\psi_j>\mu_j\},\quad \lambda_j|A_j|=1.
\]
The Robin function $R_D$ has a unique critical point $x_D$, with
$D^2R_D(x_D)>0$; this is the Robin-function input to
\cite[Theorem~1.3]{CaoGuoPengYan} and requires no strict boundary
curvature assumption.

Turkington's concentration theorem \cite{TurkingtonI}, stated in
\cite[Theorem~2.1]{CaoWangStability}, gives, after subsequence extraction,
\[
 x_j:=\int_Dx\omega_j\dd\longrightarrow x_D,\qquad
 \operatorname{diam}A_j\leq C\lambda_j^{-1/2}\longrightarrow0.
\]
Thus
\[
 A_j\subset B_{|x_j-x_D|+C\lambda_j^{-1/2}}(x_D)
 \subset B_{o(1)}(x_D).
\]
Since $H_D$ is uniformly bounded near $(x_D,x_D)$ and
$\int_D\omega_j=1$, for $x\in\partial A_j$,
\[
 \mu_j=\psi_j(x)=\int_DG_D(x,y)\omega_j(y)\,\mathrm dy
 \geq\frac{\log\lambda_j}{4\pi}-C\longrightarrow+\infty.
\]
Turkington's single-core result
\cite[Corollary~4.6, p.~1026]{TurkingtonI} gives, for large $j$,
\[
 A_j\text{ is simply connected},\qquad \partial A_j\in C^1.
\]
The nonvanishing normal derivative also follows from \cref{thm:main}.
For a fixed small ball $B_d(x_D)\subset D$,
\[
 A_j=B_d(x_D)\cap\{\psi_j>\mu_j\}\subset B_d(x_D).
\]
Thus \cite[Theorem~1.2]{CaoGuoPengYan} applies with $k=1$ and
$\kappa_1=1$, including the simply connected core required in its
published formulation. Its threshold is an unknown.

If the conclusion failed, choose distinct $A_j,B_j\in\M_{m_j}$ at
$m_j\downarrow0$. Apply the preceding results to both sequences.
They concentrate at the same nondegenerate point $x_D$, with the same
$\lambda_j$ and circulation. Local uniqueness gives
$\psi_{A_j}=\psi_{B_j}$ and $A_j=B_j$ for large $j$, a contradiction.
\end{proof}

\begin{proof}[Proof of Theorem~\ref{thm:convex-finiteness}]
Existence follows from \cref{prop:turkington}. Choose
\[
 d=\min\{m_D/2,\delta_D/2,|D|/4\}>0,\qquad I=[d,|D|-d].
\]
By \cref{an:uniform-interior,end:small,end:full},
\[
 1\leq\#\M_m\leq
 \begin{cases}
  N_I,&m\in I,\\
  1,&m\in(0,|D|)\setminus I.
 \end{cases}
\]
Hence $N_D=\max\{1,N_I\}$ depends only on $D$.
\end{proof}

\appendix
\section{General simply connected domains}\label{gen:section}
Throughout this appendix, $D$ is a bounded simply connected domain with
$C^\infty$ boundary.
The regularity, variational, and nodal arguments above apply in
any smooth bounded simply connected domain. Convexity enters through the separator in
\cref{ker:separator}, which proves $K_A^0=\{0\}$. Without it, the
bound $\dim K_A^0\leq2$ alone does not imply isolation.

A conformal change of variables makes the additional term in the
boundary comparison explicit.

Identify $\R^2$ with $\mathbb C$. Fix $A\in\M_m$ and
$0\ne\phi\in K_A^0$, and set
\[
 U=\G\one_A,\quad v=\G(\phi\,ds|_\Gamma),\quad
 P=-\partial_nU>0,\quad b=\partial_nv.
\]
Let $F:\mathbb D\to D$ be conformal, smoothly extended to the closed
disk, with $F'\ne0$ there. Put
\[
 \alpha=|F'|^2,\quad \Omega=F^{-1}(A),\quad
 \widetilde U=U\circ F,\quad \widetilde v=v\circ F.
\]
Then
\[
 -\Delta_w\widetilde U=\alpha\one_\Omega,\qquad
 m=\int_\Omega\alpha\,dw.
\]
On $w=e^{i\theta}$, write
$\widetilde P=-\partial_r\widetilde U$ and
$\widetilde b=\partial_r\widetilde v$. For
$Y(w)=a+\sigma w$, $a\in\mathbb C$, $\sigma\in\R$, and
$Z(F(w))=F'(w)Y(w)$,
\[
 \begin{gathered}
 P\circ F=\widetilde P/|F'|,\qquad
 b\circ F=\widetilde b/|F'|,\qquad
 ds_{\partial D}=|F'|\,d\theta,\\
 (Z\cdot n)\circ F=|F'|y,\qquad
 y(\theta)=\operatorname{Re}(ae^{-i\theta})+\sigma.
 \end{gathered}
\]
By \cref{ker:boundaryzeros}, $\widetilde b$ has two simple zeros.
Choose lifts $\theta_1<\theta_2<\theta_1+2\pi$ with
$\widetilde b>0$ on $(\theta_1,\theta_2)$, and take
\[
 c=\frac{\theta_1+\theta_2}{2},\quad
 \beta=\frac{\theta_2-\theta_1}{2}\in(0,\pi),\quad
 Y(w)=e^{ic}-\cos\beta\,w.
\]
Thus $y=\cos(\theta-c)-\cos\beta$ has the sign of $\widetilde b$, and
\[
 \int_{\partial D}Pb(Z\cdot n)\ds
 =\int_0^{2\pi}\widetilde P\,\widetilde b\,y\,d\theta>0.
\]

The transformed divergence and zero moment are
\[
 \begin{gathered}
 (\operatorname{div}_xZ)\circ F
 =2\operatorname{Re}\left(\sigma+\frac{F''}{F'}Y\right),\qquad
 \alpha(\operatorname{div}_xZ)\circ F
 =\operatorname{div}_w(\alpha Y),\\
 \int_\Omega\alpha\widetilde v\,dw
 =\int_Av\dd=\int_\Gamma U\phi\ds=q\int_\Gamma\phi\ds=0.
 \end{gathered}
\]
Consequently \eqref{ker:bilinear} becomes
\[
 0<\int_{\partial D}Pb(Z\cdot n)\ds
 =-\int_\Omega\widetilde v\operatorname{div}_w(\alpha Y)\,dw
 =-\int_\Omega\widetilde v\,Y\cdot\nabla\alpha\,dw.
\]
The mass constraint does not cancel this weighted moment.

For conformal fields on the physical domain,
\[
 \operatorname{div}Z\equiv2\sigma
 \quad\Longleftrightarrow\quad
 Z(z)=a+(\sigma+i\omega)z,\qquad
 a\in\mathbb C,\quad\sigma,\omega\in\R,
\]
since a holomorphic derivative with constant real part is constant.
Such fields satisfy $\int_{\partial D}Pb(Z\cdot n)\ds=0$, but the
required sign comparison is established only for convex $D$.
A conformal field tangent to $\partial D$ has zero boundary and bulk
terms in \eqref{ker:bilinear}, so adding it does not alter the obstruction.

For a general smooth field, the bilinear stress identity instead gives
\[
 \begin{aligned}
 \int_{\partial D}Pb(Z\cdot n)\ds
 &=-\int_Av\operatorname{div}Z\dd
   -\int_D(\mathcal E_Z\nabla U)\cdot\nabla v\dd,\\
 \mathcal E_Z&=DZ+(DZ)^T-(\operatorname{div}Z)I.
 \end{aligned}
\]
Prescribing divergence therefore introduces a strain term whose sign is
not controlled by the present trace inequality.

For smooth periodic families, the same stress identity includes a
term arising from the variation of the comparison field.

\begin{proposition}\label[proposition]{gen:families}
Let $A_t\in\M_m$ be a $C^2$ family of period $T$ in local
$C^{1,\alpha}$ normal graph coordinates. Set
\[
 U_t=\G\one_{A_t},\quad u_t=U_t-q(t),\quad
 \Gamma_t=\partial A_t,\quad g_t=|\nabla U_t|_{\Gamma_t},
\]
and let $\phi_t$ be the outward normal velocity. Then $q\in C^2$ and
\begin{equation}\label{gen:family-derivatives}
 \begin{gathered}
 \int_{\Gamma_t}\phi_t\ds=0,\qquad
 w_t:=\dot U_t=\G(\phi_t\,ds|_{\Gamma_t}),\qquad
 v_t:=\dot u_t=w_t-\dot q(t),\\
 S_{\Gamma_t}\phi_t-g_t\phi_t=\dot q(t),\qquad
 \phi_t\in K_{A_t},\\
 \int_{A_t}w_t\dd=0,\qquad
 \int_{A_t}v_t\dd=-m\dot q(t).
 \end{gathered}
\end{equation}
For $P_t=-\partial_nU_t>0$ and $b_t=\partial_nw_t$,
$\dot P_t=-b_t$. Every $C^1$ periodic family $Z_t$ of smooth conformal
fields on $\overline D$ satisfies
\begin{equation}\label{gen:periodic}
 \int_0^T\!\int_{\partial D}P_tb_t(Z_t\cdot n)\ds\,dt
 =\int_0^T\!\int_{A_t}u_t\operatorname{div}\dot Z_t\dd\,dt.
\end{equation}
Both sides equal
\begin{equation}\label{gen:periodic-state}
 \frac12\int_0^T\!\int_{\partial D}P_t^2(\dot Z_t\cdot n)\ds\,dt.
\end{equation}
\end{proposition}

\begin{proof}
In a fixed normal chart, write
\[
 \Xi(\xi,r)=x(\xi)+r\nu_*(\xi),\qquad
 J(\xi,r)=1+\kappa(\xi)r,\qquad A_h=\{r<h(\xi)\},
\]
and define
\[
 \langle T_hk,z\rangle
 =\int_{\Gamma_*}k(\xi)J(\xi,h(\xi))
            z(\Xi(\xi,h(\xi)))\ds_*.
\]
The trace and normal-slice estimates of \cref{prop:shape} give, uniformly
on a small graph ball,
\[
 \begin{aligned}
 \|\one_{A_{h+k}}-\one_{A_h}-T_hk\|_{\dot H^{-1}_D}
 &\leq C(\|k\|_\infty^{3/2}+\|k\|_\infty^2),\\
 \|T_h-T_\ell\|_{\mathcal L(C^{1,\alpha},\dot H^{-1}_D)}
 &\leq C(\|h-\ell\|_\infty^{1/2}+\|h-\ell\|_\infty).
 \end{aligned}
\]
Thus $h\mapsto\one_{A_h}$ is $C^1$ into $\dot H^{-1}_D$. Since
\[
 j_h=(J(\xi,h(\xi))^2+|\partial_\xi h|^2)^{1/2},\qquad
 (\phi_t\circ X_h)j_h=\dot hJ_h,
\]
we obtain $T_h\dot h=\phi_t\,ds|_{\Gamma_t}$ and the formula for
$w_t$ in $H_0^1(D)$.

The material trace is analytic by \eqref{an:layer-formula}.
Its constant value is $q(t)$, so $q\in C^2$; differentiating it gives
\[
 S_{\Gamma_t}\phi_t-g_t\phi_t=\dot q(t).
\]
Equivalently, the fixed-boundary trace derivative is
$S_{\Gamma_t}\phi_t$, and the evaluation derivative is
$\nabla U_t\cdot(\phi_t\nu_t)=-g_t\phi_t$.
Therefore
\[
 Q_{A_t}(\phi_t)=-\dot q(t)\int_{\Gamma_t}\phi_t\ds=0,\qquad
 \int_{A_t}w_t\dd
 =\int_{\Gamma_t}U_t\phi_t\ds=q(t)\int_{\Gamma_t}\phi_t\ds=0,
\]
which proves \eqref{gen:family-derivatives}.

A finite chart cover of the parameter circle gives
$\bigcup_t\overline{A_t}\subset K\subset D$. In a fixed outer collar
$\mathcal O$, the Green kernel and its derivatives are smooth for
$y\in K$. Differentiating the chart volume integrals twice yields
\[
 t\mapsto U_t\in
 C^2(\R/T\mathbb Z;C^{j+1}(\overline{\mathcal O}))
 \quad(j\geq0),\qquad \dot P_t=-b_t.
\]
No second source derivative in $\dot H^{-1}_D$ is used.

For a fixed conformal field $Z$, the state stress vector
$(Z\cdot\nabla U_t)\nabla U_t-\tfrac12|\nabla U_t|^2Z$ gives
\begin{equation}\label{gen:state-pohozaev}
 \frac12\int_{\partial D}P_t^2(Z\cdot n)\ds
 =-\int_{A_t}Z\cdot\nabla u_t\dd
 =\int_{A_t}u_t\operatorname{div}Z\dd.
\end{equation}
Since $u_t|_{\Gamma_t}=0$, differentiation gives
\[
 \int_{\partial D}P_tb_t(Z\cdot n)\ds
 =-\int_{A_t}v_t\operatorname{div}Z\dd.
\]
For $J(t)=\int_{A_t}u_t\operatorname{div}Z_t\dd$,
\[
 J'(t)
 =-\int_{\partial D}P_tb_t(Z_t\cdot n)\ds
   +\int_{A_t}u_t\operatorname{div}\dot Z_t\dd,\qquad
 J(T)=J(0).
\]
Integration proves \eqref{gen:periodic}; \eqref{gen:state-pohozaev}
with $Z=\dot Z_t$ proves \eqref{gen:periodic-state}.
\end{proof}

For translations, rotations, and dilations,
\[
 Z(z)=a+(\sigma+i\omega)z
 \quad\Longrightarrow\quad
 \int_{\partial D}P_tb_t(Z\cdot n)\ds=2\sigma m\dot q(t).
\]
In particular, dilation orthogonality at a fixed parameter requires
$\dot q(t)=0$.

\begin{proposition}\label[proposition]{gen:zero-tracking}
In \cref{gen:families}, assume $q$ is constant and $\phi_t\ne0$ for
every $t$. Fix $F:\mathbb D\to D$ conformal and put
\[
 \mathcal B(t,\theta)=|F'(e^{i\theta})|b_t(F(e^{i\theta})).
\]
Its zeros admit $C^1$ periodic circle-valued labels
$z_\pm=e^{i\theta_\pm}$, distinguished by
$\pm\partial_\theta\mathcal B(t,\theta_\pm)>0$.
They are uniformly separated. Set
\[
 a=\min_{\mathcal B=0}|\partial_\theta\mathcal B|>0,\qquad
 M=\max|\partial_t\mathcal B|.
\]
Local angular lifts satisfy
\[
 \dot\theta_\pm=r_\pm:=
 -\frac{\partial_t\mathcal B(t,\theta_\pm)}
       {\partial_\theta\mathcal B(t,\theta_\pm)},\qquad
 |r_\pm|\leq M/a.
\]
With
\[
 \ell=\operatorname{Arg}_{(0,2\pi)}(z_-/z_+),\quad
 \beta=\ell/2,\quad E=z_+e^{i\beta},
\]
one has
\[
 \dot\ell=r_--r_+,\qquad
 \dot E=\tfrac i2(r_++r_-)E.
\]
The labels have the same integer winding over a period. The fields
\[
 Y_t(w)=E(t)-\cos\beta(t)w,\qquad
 Z_t(F(w))=F'(w)Y_t(w)
\]
are $C^1$ and periodic in every $C^j(\overline D)$, with
\[
 \sup_t(\|Z_t\|_{C^j}+\|\dot Z_t\|_{C^j})
 \leq C_{F,j}(1+M/a).
\]
Their boundary integrand in \eqref{gen:periodic} is positive away
from the two zeros.
\end{proposition}

\begin{proof}
By \eqref{gen:family-derivatives},
$0\ne\phi_t\in K_{A_t}^0$; hence \cref{ker:boundaryzeros} gives two
simple zeros. Joint $C^1$ regularity and compactness of the parameter
torus imply $a>0$. Uniform continuity of
$\partial_\theta\mathcal B$ gives a common separation $\rho>0$.
The implicit-function branches glue uniquely by the derivative signs,
including at $t=0,T$, and
\[
 \rho\leq\ell\leq2\pi-\rho,\qquad
 \theta_-(T)-\theta_-(0)
 =\theta_+(T)-\theta_+(0)=2\pi k,\quad k\in\mathbb Z.
\]
Differentiate $\theta_-=\theta_++\ell$ and $E=z_+e^{i\ell/2}$
to obtain the formulas. Smoothness of $F,F^{-1}$ gives the field bounds.
Finally, with $E=e^{ic}$ locally,
\[
 Y_t\cdot n=\cos(\theta-c)-\cos\beta
 \quad\text{has the sign of }\mathcal B(t,\theta).
\]
\end{proof}

The resulting positive period integral equals the moving-field term
\eqref{gen:periodic-state}; periodicity does not force either to vanish.
No existence of such a maximizing family is asserted.

The analytic reduction also applies without convexity. Its
finite-dimensional kernel describes the remaining obstruction to
finiteness.

Compactness, graph convergence, and shape analyticity in
\cref{an:compactness,an:analytic-map} require no convexity.
For $A\in\M_m$, set
\[
 B=C^{1,\alpha}(\Gamma),\quad \mathcal X=\mathcal Y=B\times\R,\quad
 \widehat{\mathscr F}(h,\eta)=\mathscr F(h,q_A+\eta,m).
\]
Then
\[
 \widehat L(\phi,\eta)
 =\left(S_\Gamma\phi-g\phi-\eta,\int_\Gamma\phi\ds\right).
\]
For $\widehat L_0(\phi,\eta)=(-g\phi-\eta,\int_\Gamma\phi\ds)$,
\[
 \widehat L_0^{-1}(f,b)=(\phi,\eta),\qquad
 \eta=-\frac{b+\int_\Gamma f/g\ds}{\int_\Gamma1/g\ds},
 \quad \phi=-\frac{f+\eta}{g}.
\]
Since $S_\Gamma:B\to B$ is compact, $\widehat L$ is Fredholm of
index zero. Polarization and \cref{ker:bootstrap} give
\[
 \ker\widehat L
 =\{(\phi,c_A(\phi)):\phi\in K_A\},\qquad
 d:=\dim\ker\widehat L=\dim K_A
 \leq\dim K_A^0+1\leq3.
\]
Choose closed splittings
\[
 \mathcal X=N\oplus\mathcal X_1,\qquad
 \mathcal Y=\operatorname{Ran}\widehat L\oplus\mathcal C,\qquad
 N=\ker\widehat L,\quad\dim\mathcal C=d,
\]
and let $\Pi$ project onto $\operatorname{Ran}\widehat L$.
The analytic implicit-function theorem gives a local analytic map
$\Psi:N\to\mathcal X_1$, $\Psi(0)=D\Psi(0)=0$, such that
\[
 \begin{gathered}
 R_{\rm red}(k)
 =(I-\Pi)\widehat{\mathscr F}(k+\Psi(k)),\\
 \widehat{\mathscr F}(k+x)=0
 \quad\Longleftrightarrow\quad
 x=\Psi(k),\quad R_{\rm red}(k)=0.
 \end{gathered}
\]
Thus the stationary set is locally the analytic zero set of
$R_{\rm red}:\R^d\to\R^d$, $d\leq3$; it need not be a manifold.

\begin{proposition}\label[proposition]{gen:conditional}
Let $D$ be smooth, bounded, and simply connected.
If $K_A^0=\{0\}$ for every $A\in\M_m$ at a fixed mass, then
$\#\M_m<\infty$. If this exclusion holds at every maximizing pair
with mass in a compact interval $I\subset(0,|D|)$, then
$\sup_{m\in I}\#\M_m<\infty$.
If it holds at every admissible mass and the Robin function has finitely
many nondegenerate global minima, then
$\sup_{0<m<|D|}\#\M_m<\infty$.
\end{proposition}

\begin{proof}
The condition $K_A^0=0$ makes the shape-only derivative injective,
with index $-1$, so \cref{an:scalar} applies. Were $\M_m$ infinite,
its compact derived set $\M_m'$ would contain a threshold maximum
$A_0$. Distinct maximizers approaching $A_0$ give distinct zeros
$q_j\to q_0$ of the scalar residual. Therefore
\[
 r(\,\cdot\,,m)\equiv0\text{ near }q_0
 \quad\Longrightarrow\quad
 A_{H(t,m)}\in\M_m',\quad q_{A_{H(t,m)}}=t
\]
by \cref{an:constant-energy}. Taking $t>q_0$ is a contradiction.
This proves fixed-mass finiteness as in \cref{an:finiteness}.

The residual consequently has finite vanishing order at every
maximizing pair. Weierstrass preparation and joint compactness give
the bound on $I$ by the proof of \cref{an:uniform-interior}.
Near full mass, \cref{end:full} gives uniqueness without convexity.
If the global Robin minima are $p_1,\ldots,p_J$, Turkington's concentration
theorem and single-core conclusion
\cite[Corollary~4.6]{TurkingtonI}, applied as in \cref{end:small}, place
every small-mass maximizer in one of their single-patch regimes.
Nondegeneracy and local uniqueness
\cite{CaoGuoPengYan} give $\#\M_m\leq J$ for all sufficiently small
$m$. These endpoint bounds and an interior interval prove the last claim.
\end{proof}

The same local argument shows
\[
 \#\M_m=\infty,\quad
 A_0\in\operatorname*{arg\,max}_{A\in\M_m'}q_A
 \quad\Longrightarrow\quad K_{A_0}^0\ne0.
\]
For a $C^1$ fixed-mass stationary arc through a maximizer,
\[
 \frac{d}{dt}\cE(\one_{A_t})
 =q(t)\frac{d}{dt}|A_t|=0,
\]
so the entire connected arc consists of maximizers. Its velocity satisfies
\[
 \phi_t\in K_{A_t},\qquad
 \dot q(t)=0,\quad\phi_t\ne0
 \quad\Longrightarrow\quad 0\ne\phi_t\in K_{A_t}^0.
\]
For an immersed maximizing manifold $\mathcal N$ and a critical point
$A$ of $q|_{\mathcal N}$, normal velocities give an injection
\[
 T_A\mathcal N\hookrightarrow K_A^0,\qquad
 \dim\mathcal N\leq2.
\]
These are conditional and first-order statements; unconditional
finiteness for general simply connected domains is not established.

\par\noindent\begin{minipage}{\textwidth}

\section*{Statements and Declarations}

\noindent\textbf{Funding.}
\textit{ D. Cao was supported by National Key R\&D Program (Grant 2023YFA1010001) and NNSF of China (Grant
12371212). G. Wang was supported by NNSF of China (Grant 12471101) and  Fundamental Research Funds for the Central Universities (Grant DUT26YC106). Y. Wang was supported by NNSF of China (Grant 12471110). 
}


\medskip\noindent\textbf{Data availability.}
No datasets were generated or analysed during this study.


\medskip\noindent\textbf{Use of AI-assisted tools.}
During preparation of this manuscript, ChatGPT 6.0  (OpenAI) was used for
assistance with structural and language revision, LaTeX editing, and checking parts of the mathematical exposition. The authors are responsible for verifying all arguments and references and for the final content of the manuscript.

\end{minipage}

\end{document}